\documentclass[10pt]{article}

\RequirePackage{amsmath}
\RequirePackage{amsfonts}
\RequirePackage{amssymb}
\RequirePackage{hyperref}
\RequirePackage{multirow}
\RequirePackage{algorithm}
\RequirePackage{algpseudocode}
\RequirePackage{graphicx}
\RequirePackage{xcolor}
\RequirePackage{amsmath}
\RequirePackage{amsfonts}
\RequirePackage{amssymb}
\RequirePackage{hyperref}
\RequirePackage{multirow}
\RequirePackage{algorithm}
\RequirePackage{algpseudocode}
\RequirePackage{graphicx}
\RequirePackage{xcolor}
\allowdisplaybreaks

\usepackage{dpr_fec,ifthen,dsfont,tikz}
\newcommand{\ones}{\mathds{1}}
\newcommand{\newlinefortechreport}{\\}

\usepackage{isetechreport}
\def\reportyear{25T}
\def\reportno{005}
\def\revisionno{0}
\def\originaldate{\today}
\def\revisiondate{\today}
\coraltrue
\cvcrfalse

\begin{document}

\title{NonOpt: Nonconvex, Nonsmooth Optimizer\thanks{Supported by the U.S.~National Science Foundation under award number CCF-2139735 and by the Office of Naval Research under award number N00014-24-1-2703.}}

\author{Frank E.~Curtis\thanks{E-mail: frank.e.curtis@lehigh.edu}}
\author{Lara Zebiane\thanks{E-mail: lara.zebiane@lehigh.edu}}
\affil{Department of Industrial and Systems Engineering, Lehigh University}
\titlepage

\maketitle

\begin{abstract}
  NonOpt, a C++ software package for minimizing locally Lipschitz objective functions, is presented.  The software is intended primarily for minimizing objective functions that are nonconvex and/or nonsmooth.  The package has implementations of two main algorithmic strategies: a gradient-sampling and a proximal-bundle method.  Each algorithmic strategy can employ quasi-Newton techniques for accelerating convergence in practice.  The main computational cost in each iteration is solving a subproblem with a quadratic objective function, a linear equality constraint, and bound constraints.  The software contains dual active-set and interior-point subproblem solvers that are designed specifically for solving these subproblems efficiently.  The results of numerical experiments with various test problems are provided to demonstrate the speed and reliability of the software.
\end{abstract}

\renewcommand{\labelitemi}{$\circ$}
\renewcommand{\labelitemii}{$\bullet$}

\section{Introduction}\label{sec.introduction}

This paper provides an introduction to an open-source C++ software package known as NonOpt, which stands for \emph{Nonconvex, Nonsmooth Optimizer} and is available at \href{https://github.com/frankecurtis/NonOpt}{\texttt{https://github.com/frankecurtis/NonOpt}}.  This software has implementations of a suite of algorithms for solving continuous optimization problems when the objective function is nonconvex and/or nonsmooth.  The implemented algorithms are applicable for solving unconstrained problems, which means they could also be employed effectively as subproblem solvers within various algorithmic approaches for solving other types of optimization problems, say, involving discrete decision variables or constraints.  The main algorithmic strategies that are implemented in the NonOpt software are a gradient-sampling method and a proximal-bundle method.  The software comes equipped with options for employing quasi-Newton techniques---specifically, Davidon-Fletcher-Powell (DFP) \cite{Davi91} and Broyden-Fletcher-Goldfarb-Shanno (BFGS) techniques \cite{Broy70,Flet70,Gold70,Shan70}---within the algorithms.  Convergence guarantees can be maintained with such techniques, and importantly their use significantly improves practical performance.

NonOpt builds from earlier Matlab implementations of algorithms that were written and employed for the numerical experiments in \cite{CurtOver12,CurtQue13,CurtQue15}.  An early C++ implementation of the software (known as SVANO) was written and employed for the experiments in \cite{CurtRobiZhou20}.  Further improvements to the software led to the implementation (referred to as NonOpt) written and employed for the experiments in~\cite{CurtLi22}.  The present paper provides a thorough description of an updated version of the NonOpt software (version 2.0), which has not appeared previously.  It also highlights a significant new feature to the software, namely, an interior-point method based on the strategy in \cite{CurtNoce07} for solving the quadratic optimization (QP) subproblems that arise in gradient-sampling and proximal-bundle methods.  The overall algorithmic strategies can make use of inexact subproblem solutions obtained from these solvers by following the procedures in \cite{CurtLi22}, which is particularly effective when solving large-scale problems and the new interior-point subproblem solver is employed.  This new interior-point subproblem solver complements a previously implemented dual active-set method.  Combined, with these two strategies for solving the arising QP subproblems, the software can solve both small- and large-scale instances efficiently.  The numerical experiments that are presented in this paper demonstrate the practical performance of the gradient-sampling and proximal-bundle methods that are implemented in the NonOpt software.

An important feature of the open-source NonOpt software is that it has been written to be extensible, allowing a user to include other search-direction computation schemes, globalization mechanisms, or (quasi-)Newton-based strategies.

\subsection{Organization}

A formal statement of the optimization problems that the algorithms implemented in NonOpt are designed to solve, as well as preliminary mathematical background on stationarity conditions and theoretical guarantees for the implemented algorithms, are provided in \S\ref{sec.problem}.  A high-level description of the implemented algorithmic framework, which includes both gradient-sampling and proximal-bundle strategies as special cases, is presented in \S\ref{sec.algorithm}.  This section also provides a detailed account of the modules in the code through which one distinguishes between different search-direction computation strategies (e.g., gradient-sampling versus proximal-bundle), globalization mechanisms, quasi-Newton strategies, and others.  A complete description of the newly implemented interior-point method for solving the arising QP subproblems is presented in \S\ref{sec.subproblem}.  This method has been tailored specifically for the form of the QPs that arise in NonOpt algorithms.  Further details of the software implementation are presented in \S\ref{sec.implementation} and the results of extensive numerical experiments are presented in \S\ref{sec.numerical}.  Concluding remarks are provided in \S\ref{sec.conclusion}.

\section{Problem Formulation and Preliminaries}\label{sec.problem}

NonOpt is designed to solve optimization problems of the form
\bequation\label{prob.opt}
  \min_{x \in \R{n}} f(x),
\eequation
where the objective function $f$ satisfies the following assumption.

\bassumption\label{ass.basic}
  The objective function $f : \R{n} \to (\R{} \cup \{\infty\})$ has an effective domain
  \bequationNN
    \dom(f) := \{x \in \R{n} : f(x) < \infty\}
  \eequationNN
  that is nonempty with a nonempty interior, i.e., $\interior(\dom(f)) \neq \emptyset$.  In addition, $f$ is locally Lipschitz over $\interior(\dom(f))$ in the sense that for every $x \in \interior(\dom(f))$ there exist positive real numbers $\epsilon_{f,x} \in \R{}_{>0}$ and $L_{f,x} \in \R{}_{>0}$ such that, for all $(\xbar,\xhat) \in \R{n} \times \R{n}$ with $\|\xbar - x\|_2 \leq \epsilon_{f,x}$ and $\|\xhat - x\|_2 \leq \epsilon_{f,x}$, one has that
  \bequationNN
    |f(\xbar) - f(\xhat)| \leq L_{f,x} \|\xbar - \xhat\|_2.
  \eequationNN
\eassumption

It is well known that, by Rademacher's theorem \cite{Rade1919}, the objective function $f$ being locally Lipschitz over $\interior(\dom(f))$ means that it is continuous over $\interior(\dom(f))$ and differentiable almost everywhere in $\interior(\dom(f))$ in the sense that the subset of points in $\interior(\dom(f))$ at which it is not differentiable has a Lebesgue measure of zero \cite{Clar83}.  Under Assumption~\ref{ass.basic}, the Clarke generalized directional derivative function with respect to $f$, defined for all $x \in \R{n}$ and $d \in \R{n}$ by
\bequationNN
  f^\circ(x,d) = \limsup_{\substack{\xbar \to x \\ \alpha \searrow 0}} \frac{f(\xbar + \alpha d) - f(\xbar)}{\alpha},
\eequationNN
defines for any $x \in \interior(\dom(f))$ a mapping $f^\circ(x,\cdot) : \R{n} \to \R{}$ that is finite-valued and sublinear with $|f^\circ(x,d)| \leq L_{f,x} \|d\|_2$ for all $d \in \R{n}$, where $L_{f,x} \in \R{}_{>0}$ is a Lipschitz constant for $f$ near $x$.  In addition, for any $x \in \R{n}$, the set of Clarke generalized gradients for the objective function $f$ at $x$ is given by
\begin{align*}
  \gengradient f(x) :=&\ \{g \in \R{n} : f^\circ(x,d) \geq g^Td\ \text{for all}\ d \in \R{n}\} \\
  =&\ \conv(\{ g \in \R{n} : \{\nabla f(x_k)\} \to g\ \text{for some}\ \{x_k\} \to x\ \text{with}\ \{x_k\} \subset \Dcal_{\nabla f} \}),
\end{align*}
where $\Dcal_{\nabla f} \subseteq \interior(\dom(f))$ is the full-measure subset over which $f$ is differentiable.  The latter equation, which states that the Clarke generalized gradient of $f$ at $x$ is equal to the convex hull of limits of gradient sequences, is well known; see \cite{Clar83,RockWets1998}.

The algorithms that are implemented in NonOpt are well defined under the relatively loose Assumption~\ref{ass.basic} as long as an initial iterate $x_1$ is provided that lies in $\interior(\dom(f))$.  However, it should be mentioned that to attain the theoretical convergence guarantees offered for NonOpt's implemented gradient-sampling or proximal-bundle method, the objective function should satisfy additional assumptions as well.  First, the objective function should in fact be real-valued over $\R{n}$.  NonOpt has safeguards to ensure that no iterate (after $x_1$) is generated at which the objective function is infinite-valued, although technically speaking the implemented algorithms might not be guaranteed to have convergence guarantees when~$f$ can take infinite values.  Second, the objective function should satisfy at least one of two additional assumptions depending on whether the gradient-sampling or proximal-bundle method is employed.  Specifically, to be consistent with the convergence guarantees offered in, e.g., \cite{BurkCurtLewiOverSimo20,BurkLewiOver2005,Kiwi2007} and/or \cite{Mifflin1982,SchrZowe92}, problem~\eqref{prob.opt} should satisfy the following assumption in addition to Assumption~\ref{ass.basic}.

\bassumption\label{ass.theory}
  The (real-valued) objective function $f : \R{n} \to \R{}$ satisfies at least one of the following two conditions:
  \benumerate
    \item[(a)] $f$ is continuously differentiable on an open set $\Dcal_{\nabla f}$ with full measure over $\R{n}$;
    \item[(b)] $f$ is weakly upper semismooth over $\R{n}$ in the sense that, for all $x \in \R{n}$ and $d \in \R{n}$, the directional derivative of $f$ at $x$ along $d$, namely,
    \bequationNN
      f'(x,d) := \lim_{\alpha \searrow 0} \frac{f(x + \alpha d) - f(x)}{\alpha}
    \eequationNN
    exists and for all sequences $\{(\alpha_k,g_k)\} \subset \R{}_{>0} \times \R{n}$ satisfying $\{\alpha_k\} \searrow 0$ and $g_k \in \gengradient f(x + \alpha_k d)$ for all $k \in \N{}$ the directional derivative satisfies
    \bequationNN
      \liminf_{k \to \infty} g_k^Td \geq f'(x,d).
    \eequationNN
  \eenumerate
\eassumption

\noindent
Assumption~\ref{ass.theory}(a) ensures that convergence guarantees for the implemented gradient-sampling methods hold.  On the other hand, Assumption~\ref{ass.theory}(b) ensures guarantees for the implemented proximal-bundle methods.

Given subroutines for evaluating the objective~$f$ and an element of $\gengradient f(x)$ for any $x \in \interior(\dom(f))$ in addition to an initial iterate $x_1 \in \interior(\dom(f))$, NonOpt generates an iterate sequence $\{x_k\} \subset \dom(f)$ that is intended to converge to stationarity with respect to $f$.  That is, at the very least, the software aims to generate such a sequence $\{x_k\}$ with a limit point $\xbar$ satisfying $0 \in \gengradient f(\xbar)$.

\section{Algorithmic Framework}\label{sec.algorithm}

In this section we provide a high-level description of the algorithmic framework for solving instances of problem~\eqref{prob.opt} that is implemented in NonOpt.  The framework encapsulates various instances of gradient-sampling and proximal-bundle methods, where to employ each particular instance a user can select between multiple options related to the line search, quasi-Newton strategy, and other aspects.  We begin this section with a mathematical description of the main components of the framework, then discuss the modules of the code that allow a user to switch between different algorithmic strategies within the gradient-sampling and proximal-bundle methods.

\subsection{Main algorithmic components}

At initialization with a user-provided initial point $x_1 \in \R{n}$, the objective-function value $f_1 := f(x_1)$, a generalized gradient $g_1 \in \gengradient f(x_1)$, and a symmetric positive definite matrix $H_1$ with inverse $W_1 := H_1^{-1}$ are all used to compute a search direction.  More generally, at the beginning of iteration $k \in \N{}$ and with the notation $[m] := \{1,\dots,m\}$ for any $m \in \N{}$, the algorithm has accumulated and/or will generate a set of points $\{x_{k,j}\}_{j\in[m]} \subset \R{n}$, real numbers $\{f_{k,j}\}_{j\in[m]} \subset \R{}$, and real vectors $\{g_{k,j}\}_{j\in[m]} \subset \R{n}$ that are used along with symmetric positive definite $H_k$ and/or $W_k := H_k^{-1}$ to compute the search direction $d_k \in \R{n}$.  In a gradient-sampling method, the set of points $\{x_{k,j}\}_{j\in[m]} \subset \R{n}$ corresponds to the current iterate~$x_k$, previous iterates and sample points that are contained within a neighborhood of the current iterate, and newly generated sample points within the same neighborhood; the set of real numbers $\{f_{k,j}\}_{j\in[m]} \subset \R{n}$ has each element equal to the current objective function value $f(x_k)$; and the real vectors $\{g_{k,j}\}_{j\in[m]} \subset \R{n}$ correspond to gradients of~$f$ at the points in $\{x_{k,j}\}_{j\in[m]}$.  (The strategies for maintaining and generating these points follow the strategies proposed for gradient-sampling methods in \cite{BurkCurtLewiOverSimo20,BurkLewiOver2005,CurtQue13,CurtQue15,CurtRobiZhou20,Kiwi2007}.)  In a proximal-bundle method, the set of points $\{x_{k,j}\}_{j\in[m]}$ contains the current iterate~$x_k$ and (a subset of) prior iterates, the set of real numbers $\{f_{k,j}\}_{j\in[m]} \subset \R{}$ correspond to the objective-function values at these points (potentially modified through \emph{downshifting} with the $\alpha$-function for an objective function stated on page~81 of~\cite{Mifflin1982}), and the set of real vectors correspond to generalized gradients of $f$ at the points $\{x_{k,j}\}_{j\in[m]}$.  (The strategy for maintaining and generating these points follows \cite{Mifflin1982,SchrZowe92}.)

Regardless of the source of $\{(x_{k,j},f_{k,j},g_{k,j})\}_{j\in[m]}$, the search direction computation employs the piecewise linear model $l_{k,m} : \R{n} \to \R{}$ defined by
\bequationNN
  l_{k,m}(x) = \max_{j \in [m]} \{f_{k,j} + g_{k,j}^T(x - x_{k,j})\}.
\eequationNN
Given $H_k \succ 0$, a corresponding piecewise quadratic model is $q_{k,m} : \R{n} \to \R{}$ with
\bequationNN
  q_{k,m}(x) = l_{k,m}(x) + \thalf (x - x_k)^TH_k(x - x_k).
\eequationNN
A step toward minimizing $f$ is given by the minimizer of $q_{k,m}$ within a trust region with norm $\|\cdot\|_{\infty}$ and radius $\delta_k \in \R{}_{>0}$, which leads to the subproblem
\bequationNN
  \min_{x \in \R{n}}\ q_{k,m}(x)\ \ \st\ \ x \in \Xcal_k := \{x \in \R{n} : \|x - x_k\|_\infty \leq \delta_k\}.
\eequationNN
This subproblem has a convex feasible region and a strongly convex objective function, so it has a unique optimal solution.  (The trust-region radius is always chosen to be finite, although depending on the choices of a couple of critical input parameters, it is unlikely to be active at the optimal solution; see \S\ref{sec.implementation}.)  With the change of variables $d \equiv x - x_k$, an equivalent smooth subproblem is
\bequation\label{eq.primal_subproblem}
  \baligned
    \min_{(d,z) \in \R{n} \times \R{}} &\ z + \thalf d^TH_kd \\
    \st &\ x_k + d \in \Xcal_k\ \ \text{and}\ \ f_{k,j} + g_{k,j}^T(x_k + d - x_{k,j}) \leq z\ \ \text{for all}\ \ j \in [m].
  \ealigned
\eequation
Let us call this the primal subproblem.  Rather than solve it directly, NonOpt has implemented solvers (see~\S\ref{sec.subproblem}) for the corresponding dual subproblem given by
\bequation\label{eq.dual_subproblem}
  \begin{aligned}
    \max_{(\omega,\gamma) \in \R{m} \times \R{n}} &\ -\thalf (G_k\omega + \gamma)^TW_k(G_k\omega + \gamma) + b_k^T\omega - \delta_k \|\gamma\|_1 \\
    \st &\ \ones^T\omega = 1\ \ \text{and}\ \ \omega \geq 0,
  \end{aligned}
\eequation
where $G_k := \bbmatrix g_{k,1} & \cdots & g_{k,m} \ebmatrix$, $\ones$ (here and throughout the paper) is a vector of ones of appropriate length, and, for all $j \in [m]$, the $j$th element of $b_k$ is
\bequation\label{eq.b}
  b_{k,j} = f_{k,j} + g_{k,j}^T(x_k - x_{k,j});
\eequation
see~\cite{CurtRobiZhou20}.  There is not necessarily a unique solution of \eqref{eq.dual_subproblem} due, e.g., to the fact that the columns of $G_k$ might not be linearly independent.  That being said, any optimal solution $(\omega_k,\gamma_k)$ of \eqref{eq.dual_subproblem} yields the same value of $d_k \gets -W_k(G_k\omega_k + \gamma_k)$, which is part of the unique solution $(d_k,z_k)$ of the primal subproblem.

Upon computation of the search direction $d_k$, the NonOpt framework conducts a line search to compute the subsequent iterate $x_{k+1}$.  The line search at least imposes on the computed step size $\alpha_k \in \R{}_{>0}$ the sufficient decrease condition
\bequation\label{eq.sufficient_reduction}
  f(x_{k+1}) - f(x_k) \leq - \thalf c \alpha_k (G_k\omega_k + \gamma_k)^TW_k(G_k\omega_k + \gamma_k) = -\thalf c \alpha_k \|d_k\|_{H_k}^2
\eequation
for some $c \in (0,1)$.  There is also an option for asking the line-search mechanism to attempt to satisfy a curvature condition as well, as in the Armijo-Wolfe-type line-search schemes proposed in \cite{LewiOver2013,Mifflin1982}.  In a gradient-sampling method, in order to ensure theoretical convergence guarantees, the line search technically needs to require that $f$ is differentiable at $x_{k+1}$.  However, this condition is ignored in the software so that the user does not need to provide a mechanism for checking differentiability of $f$.  It is assumed that, at such a point of nondifferentiability, a request for a gradient of $f$ returns a generalized gradient of $f$ at the point.

The other main components of the framework in each iteration $k \in \N{}$ are to select the \emph{point set} that defines $\{(x_{k+1,j},f_{k+1,j},g_{k+1,j})\}$ and to update a quasi-Newton matrix.  The NonOpt software is set up to allow other mechanisms for selecting the point set, although at this time the only mechanism that is implemented is to remove points that either (a) fall outside of a ball centered at $x_{k+1}$ or (b) are the ``oldest'' when the cardinality of the point set would otherwise exceed a prescribed limit.  The latter possibility is straightforward to implement, so we do not provide the details here.  As for (a), it is also straightforward: Defining
\bequation\label{eq.ball}
  \mathbb{B}_{\leq \epsilon_k}^n(x_k) := \{x \in \R{n} : \|x - x_k\|_2 \leq \epsilon_k\},
\eequation
the strategy removes an element $(\xbar,\fbar,\gbar)$ from the point set if $\xbar \notin \mathbb{B}_{\leq \epsilon_{k+1}}^n(x_{k+1})$.  As for the quasi-Newton updating mechanism, we follow the strategy proposed in~\cite{CurtRobiZhou20}.  For concreteness, we discuss here the strategy in terms of BFGS updating; the strategy for DFP updating is similar \cite{NoceWrig06}.  First, in order to ensure the self-correcting properties of BFGS updating, we employ the \emph{damping} strategy from~\cite{CurtRobiZhou20}, which involves replacing $y_{k+1} \gets g_{k+1} - g_k$ by a convex combination, call it $v_k$, of $s_k \gets \alpha_k d_k$ and $y_k$ such that the pair of inequalities
\bequation\label{eq.bounds}
  \eta \leq \frac{s_k^Tv_k}{\|s_k\|_2^2}\ \ \text{and}\ \ \frac{\|v_k\|_2^2}{s_k^Tv_k} \leq \psi
\eequation
holds for user-defined $\eta \in \R{}_{>0}$ and $\psi \in [\eta,\infty)$.  Then, depending on whether a \emph{full-memory} or \emph{limited-memory} \cite{Noce1980} option is chosen, the algorithm either sets
\bsubequations\label{eq.HW_update}
  \begin{align}
    H_{k+1} &\gets \(I - \frac{s_ks_k^TH_k}{s_k^TH_ks_k}\)^T H_k \(I - \frac{s_ks_k^TH_k}{s_k^TH_ks_k}\) + \frac{v_kv_k^T}{s_k^Tv_k} \\ \text{and}\ \ 
    W_{k+1} &\gets \(I - \frac{v_ks_k^T}{s_k^Tv_k}\)^TW_k\(I - \frac{v_ks_k^T}{s_k^Tv_k}\) + \frac{s_ks_k^T}{s_k^Tv_k}
  \end{align}
\esubequations
or updates a set of pairs $\{(s_j,v_j)\}$ that is used to define $H_{k+1}$ and $W_{k+1}$ implicitly in the subsequent iteration, as in a standard L-BFGS approach \cite{ByrdNoceSchn1994,LiuNoce1989,NoceWrig06}.

A statement of the framework that is implemented in NonOpt is provided as Algorithm~\ref{alg.nonopt} below.  The software involves numerous other input parameters, all of which are documented in the code itself.  The statement of Algorithm~\ref{alg.nonopt} includes only a few parameters that, in our experience, are the most important ones in terms of influencing the practical performance of the implemented algorithms.  See the next subsection and \S\ref{sec.implementation} for further details on the parameter choices.

\floatname{algorithm}{Algorithm}
\balgorithm[ht]
  \caption{NonOpt Algorithm Framework}
  \label{alg.nonopt}
  \balgorithmic[1]
    \Require $c \in (0,1)$, $\eta \in \R{}_{>0}$, and $\psi \in [\eta,\infty)$; an initial point $x_1 \in \R{n}$; radii $(\epsilon_1,\delta_1) \in \R{}_{>0} \times \R{}_{>0}$; an integer $p \in \N{}$; and symmetric positive definite $H_1 \in \R{n \times n}$ with $W_1 \gets H_1^{-1}$
    \State set $\Bcal_1 \gets \{(x_1,f(x_1),g_1)\}$ for some $g_1 \in \gengradient f(x_1)$
    \For{\textbf{all} $k \in \N{}$}
      \State (optional) for each $j \in [p]$, sample $\xbar$ from a uniform distribution over $\mathbb{B}_{\leq \epsilon_k}^n(x_k)$ and
      \bequationNN
        \Bcal_k \gets \Bcal_k \cup \{(\xbar,f(\xbar),\gbar)\}\ \ \text{for some}\ \ \gbar \in \gengradient f(\xbar)
      \eequationNN
      \State letting $\{(x_{k,j},f_{k,j},g_{k,j})\}_{j\in[m]}$, $G_k$, and $b_k$ be defined by $\Bcal_k$, solve \eqref{eq.dual_subproblem} for $(\omega_k,\gamma_k)$
      \State set $d_k \gets -W_k(G_k\omega_k + \gamma_k)$
      \State compute $\alpha_k$ to satisfy \eqref{eq.sufficient_reduction} by a backtracking Armijo or Armijo-Wolfe line search
      \State set $x_{k+1} \gets x_k + \alpha_k d_k$ and $s_k \gets \alpha_k d_k$
      \State set $\Bcal_{k+1} \gets \{(x_{k+1},f(x_{k+1}),g_{k+1})\}$ for some $g_{k+1} \in \gengradient f(x_{k+1})$
      \State set $\epsilon_{k+1} \in (0,\epsilon_k]$ and $\delta_{k+1} \in \R{}_{>0}$ \label{step.epsilon}
      \State (optional) remove $(\xbar,\fbar,\gbar)$ from $\Bcal_{k+1}$ such that $\xbar \notin \mathbb{B}_{\leq \epsilon_{k+1}}(x_{k+1})$
      \State (optional) remove from $\Bcal_{k+1}$ its ``oldest'' elements until $|\Bcal_{k+1}|$ is below a limit
      \State set $y_k \gets g_{k+1} - g_k$ and compute $\beta_k$ as the smallest value in $[0,1]$ such that \label{step.y}
      \bequationNN
        v_k \gets \beta_k s_k + (1 - \beta_k)y_k
      \eequationNN
      \State yields \eqref{eq.bounds}, then set $(H_{k+1},W_{k+1})$ by \eqref{eq.HW_update} (or update a set of pairs $\{(s_j,v_j)\}$)
    \EndFor
  \ealgorithmic
\ealgorithm

\subsection{Modules defining key algorithmic strategies}

The source code for NonOpt involves low-level and high-level components.  By low-level components, we are referring to classes for defining and manipulating objects such as vectors and symmetric matrices, which rely on BLAS \cite{BLAS1,BLAS2,BLAS3} and LAPACK \cite{lapack} routines, as well as classes for defining and manipulating options (i.e., user-modifiable input parameters) and output streams.  The low-level components also include container classes for storing iterate information and other quantities that are passed and/or referenced throughout the software.  We do not discuss these low-level components in any further detail in this paper.

The high-level components of NonOpt relate to the main algorithmic strategies.  These are managed through a \texttt{Strategies} object through which a user can select different options as building blocks for an algorithm instance.  The following is a list of the strategies in the code along with corresponding options/parameters that, in our experience, have the greatest effect on the performance of the software.
\bitemize
  \item \texttt{approximate\_hessian\_update}: Strategy for updating the sequences of Hessian approximations $\{H_k\}$ and inverse Hessian approximations $\{W_k\}$ that appear in the subproblems (see \eqref{eq.primal_subproblem}--\eqref{eq.dual_subproblem}) for computing search directions; current options are BFGS (default) and DFP.  These employ standard procedures---see, e.g., \cite{NoceWrig06}---as well as the damping strategy described in \cite{CurtRobiZhou20}.  Information is stored for maintaining both a sequence of Hessian approximations and the corresponding sequence of inverse Hessian approximations.  (The exact type of information that is stored depends on the \texttt{SymmetricMatrix} strategy, described later on.)  The most influential options related to this strategy are the following.
  \bitemize
    \item \texttt{BFGS\_correction\_threshold\_1} and \texttt{DFP\_correction\_threshold\_1} for BFGS and DFP updates, respectively.  These correspond to the parameter~$\eta$ in~\eqref{eq.bounds}.  The default value is $10^{-8}$ for both options.
    \item \texttt{BFGS\_correction\_threshold\_2} and \texttt{DFP\_correction\_threshold\_2} for BFGS and DFP updates, respectively.  These correspond to the parameter~$\psi$ in~\eqref{eq.bounds}.  The default value is $10^8$ for both options.
  \eitemize
  \item \texttt{derivative\_checker}: Strategy for checking user-supplied derivatives; current option is only to compare derivatives to finite-difference approximations.  (The user should be aware of the fact that such comparisons may be faulty when a given point is one of nondifferentiability of the objective function.)  The derivative-checker is off by default, but can be turned on for the purpose of debugging a user-provided problem instance.  The most influential option related to this strategy is the following.
  \bitemize
    \item \texttt{DEFD\_increment} is the increment or spacing used to compute the finite-difference approximation of the partial derivative of the objective function with respect to each coordinate.  The default value is $10^{-8}$.
  \eitemize
  \item \texttt{direction\_computation}:\label{part.dc} Strategy for building a quadratic optimization (QP) subproblem to be solved for computing the search direction in each iteration; current options are gradient, gradient combination (as in a gradient-sampling method), and cutting plane (default, as in a proximal-bundle method).  The latter two strategies have been mentioned previously in this section.  The former (gradient) strategy constructs the subproblem using only the current iterate, objective function value, and generalized gradient value, as in a method for smooth optimization; it offers no convergence guarantees in general.  (The resulting search direction still employs a Hessian approximation, so this gradient strategy could be used, e.g., to employ a straightforward BFGS-based algorithm.)  The most influential options related to this strategy are the following.
  \bitemize
    \item \texttt{DCGC\_try\_aggregation} and \texttt{DCCP\_try\_aggregation} for the gradient combination and cutting plane techniques, respectively.  These options determine whether a gradient aggregation strategy is employed, as explained in \cite{CurtLi22}.  The default value is \texttt{false} for both options, but setting it to \texttt{true} can lead to significantly improved performance when solving some problems.
    \item \texttt{DCGC\_try\_gradient\_step} and \texttt{DCCP\_try\_gradient\_step} for the gradient combination and cutting plane techniques, respectively.  These options determine whether a gradient-descent-type step (scaled by $W_k$, namely, $-W_kg_k$) is considered as the step before the full subproblem solution (see \eqref{eq.dual_subproblem}) is computed.  The default value is \texttt{true} for both options.
    \item \texttt{DCGC\_try\_shortened\_step} and \texttt{DCCP\_try\_shortened\_step} for the gradient combination and cutting plane techniques, respectively.  These options determine whether a shortened version of the step computed from a subproblem is considered before the data defining the subproblem is augmented to lead to an improved search direction.  The default value is \texttt{true} for both options.
  \eitemize
  \item \texttt{line\_search}: Strategy for computing a step size along each search direction; current options are an Armijo/weak-Wolfe line search (default) or backtracking Armijo line search.  By default, each subproblem is defined with a trust-region constraint, although through user-modifiable options a user can make the trust-region more likely or not to influence the search direction.  (By default, the parameters that influence the trust-region radii are such that the trust region is rarely active in our numerical experiments.)  Overall, through user-modifiable parameters, a user can invoke various combinations of line-search and trust-region strategies.  The most influential options are the following.
  \bitemize
    \item \texttt{LSWW\_stepsize\_initial}, \texttt{LSWW\_stepsize\_sufficient\_decrease\_threshold}, and \newlinefortechreport\texttt{LSWW\_stepsize\_curvature\_threshold} for the Armijo/weak-Wolfe line search.  These correspond to the initial step size, sufficient decrease parameter, and curvature parameter.  The default values are $1$, $10^{-10}$, and $0.9$, respectively.
    \item \texttt{LSB\_stepsize\_initial} and \texttt{LSB\_stepsize\_sufficient\_decrease\_threshold} for the backtracking line search scheme.  These correspond to the initial step size and sufficient decrease parameter.  The default values are $1$ and $10^{-10}$.
  \eitemize
  \item \texttt{point\_set\_update}: Strategy for updating the set of points (with corresponding function and/or generalized-gradient values) that may be used to construct subproblems in subsequent iterations; current option is only to update set based on proximity (i.e., distance) to the current iterate.  The most influential options related to this strategy are the following.
  \bitemize
    \item \texttt{PSP\_envelope\_factor} is the factor---i.e., multiple of $\epsilon_k$---used to define the neighborhood (see \eqref{eq.ball}) outside of which points are removed from the point set.  The default value is $10^2$.
    \item \texttt{PSP\_size\_factor} is the factor---i.e., multiple of $n$---used to define the cardinality limit for the point set.  If this limit is exceeded, then old points are removed.  The default value is $5 \times 10^{-2}$.
  \eitemize
  \item \texttt{qp\_solver\_small\_scale} and \texttt{qp\_solver\_large\_scale}: Strategies for solving the arising QP subproblems, at least approximately according to the termination tests derived from \cite{CurtLi22}; current options for each are a dual active-set method (preferable for small-scale problems) and an interior-point method (preferable for large-scale problems).  In each iteration, whether the small-scale or large-scale QP solver is employed is determined by the direction computation strategy.  (The gradient direction computation strategy always employs the small-scale solver.  The gradient combination and cutting plane strategies choose the solver in iteration $k$ based on the number of columns of the matrix $G_k$.)  The newly implemented interior-point solver is presented in detail in the next section.  The most influential options related to this strategy are discussed in \S\ref{sec.implementation}.
  \item \texttt{symmetric\_matrix}:\label{part.sm} Strategy for storing the aforementioned (inverse) Hessian approximation; current options are as a dense matrix (default, based on all prior pairs for quasi-Newton updating) or as a limited-memory matrix (involving only a prescribed number of recent pairs for quasi-Newton updating).  The most influential option related to this strategy is the following.
  \bitemize
    \item \texttt{SMLM\_history} for the limited-memory technique is the typical history length for limited-memory quasi-Newton techniques.  The default value is 20.
  \eitemize
  \item \texttt{termination}: Strategy for terminating the overall algorithm that aims to solve problem~\eqref{prob.opt}; current options are a basic strategy (default) based on the norm of the most recently computed search direction or a more expensive strategy that involves computing the minimum-norm element of the convex hull of stored gradients evaluated at points near the current iterate.  The most influential options related to this strategy are the following.
  \bitemize
    \item \texttt{TB\_objective\_similarity\_tolerance} for the basic termination strategy (and \newlinefortechreport\texttt{TS\_objective\_similarity\_tolerance} for the other strategy) determines a relative tolerance for changes in the objective function value.  If the objective function value does not change more than this threshold for a sequence of consecutive iterations, then the sampling radius is reduced or the algorithm terminates due to lack of further progress.  The default value is $10^{-5}$.
  \eitemize
\eitemize

\section{Subproblem Solvers}\label{sec.subproblem}

Solving subproblem~\eqref{eq.dual_subproblem} is the most computationally expensive aspect of each iteration of Algorithm~\ref{alg.nonopt}, at least without considering the computational expense of each (user-defined) function and generalized gradient evaluation.  For solving this subproblem in each iteration, a dual active-set algorithm was implemented previously within NonOpt (version 1.0).  That dual active-set algorithm is based on the one proposed in \cite{Kiwi1986}, although a major difference is that the method in NonOpt employs a more general distance metric, namely, that defined by the matrix $W_k$ in \eqref{eq.dual_subproblem}.  (The method proposed in \cite{Kiwi1986} only considers $W_k = I$.)  The implementation of this dual active-set algorithm was used for the experiments in~\cite{CurtLi22,CurtRobiZhou20}.

Like other active-set algorithms for solving quadratic optimization problems, a drawback of this dual active-set approach is that it may require a large number of iterations to reach an exact solution or even an inexact solution that is close to the exact solution in some sense.  A recent major addition to the NonOpt software is an interior-point algorithm for solving the arising subproblems.  As shown in our numerical experiments in \S\ref{sec.numerical}, the interior-point algorithm reduces CPU time significantly compared to the dual active-set method when solving large-scale problems.

We now provide a detailed description of the interior-point algorithm implemented in NonOpt for solving subproblems of the form~\eqref{eq.dual_subproblem}.  This subproblem solver is based on the approach proposed in~\cite{CurtNoce07} along with a Mehrotra predictor-corrector strategy~\cite{Mehr1992}.  It has also been tailored explicity for solving~\eqref{eq.dual_subproblem}.  First, due to the nonsmooth $\ell_1$-norm term in the objective function of~\eqref{eq.dual_subproblem}, we split the variable $\gamma$ into a pair of nonnegative variables that correspond to the positive and negative parts of $\gamma$, call them $\sigma$ and $\rho$, respectively.  The solution of the resulting subproblem can be obtained by solving the minimization problem
\bequation\label{eq.quad_prob}
  \min_{\theta\in\R{\ell}}\ \thalf \theta^T Q \theta + c^T\theta\ \st\ A \theta = b\ \text{and}\ \theta \geq 0,
\eequation
where $\ell := m + 2n$ and the problem variables and data are given respectively by
\bequationNN
  \theta \equiv \bbmatrix \omega \\ \sigma \\ \rho \ebmatrix,\ Q = \bbmatrix G_k^TW_kG_k & G_k^TW_k & -G_k^TW_k \\ W_k G_k & W_k & -W_k \\ -W_k G_k & -W_k & W_k \ebmatrix,\ c = \bbmatrix -b_k \\ \delta_k \ones \\\delta_k \ones \ebmatrix,\ A = \bbmatrix \ones \\ 0 \\ 0 \ebmatrix^T\!\!\!\!,\ \text{and}\ b = 1.
\eequationNN
Here, for the sake of notational simplicity, we drop the subscript $k$ in \eqref{eq.quad_prob}.  For our subsequent discussions in this section, the iterates of our implemented interior-point method for solving~\eqref{eq.quad_prob} will use iteration index $j$, which the reader should not confuse with the meaning of the index $j$ in, e.g., \eqref{eq.b}.  Let us also note that if the trust region is inactive in the optimal solution of \eqref{eq.primal_subproblem}, then the solution of~\eqref{eq.dual_subproblem} is guaranteed to have $\gamma_k = 0$, in which case the solution of~\eqref{eq.quad_prob} is guaranteed to have $\rho = \sigma = 0$.  Thus, in this case, the variables $(\rho,\sigma)$ do not even need to be introduced and~\eqref{eq.quad_prob} only needs to be defined and solved over the variable~$\omega$.  The remainder of our description of our interior-point method is effectively the same in either case, so for the sake of generality we only present our method for the case when the trust region constraint might be active.  In \S\ref{sec.implementation}, we discuss our approach for trying to avoid the introduction of these additional variables in practice.

Problem~\eqref{eq.quad_prob} is a convex QP, although as previously mentioned it might not have a unique solution.  Indeed, the matrix $Q$ is not positive definite.  (One straightforward setting in which the solution is not unique is when the matrix $G_k$ contains two columns that are equal to each other and the corresponding components of $b_k$ are also equal.)  The KKT conditions for problem~\eqref{eq.quad_prob} are
\bsubequations\label{eq.kkt}
  \begin{align}
  A\theta &= b, \label{eq.kkt1} \\
  A^T u + v -Q\theta-c &= 0, \label{eq.kkt2} \\
  \Theta v &= 0, \label{eq.kkt3} \\ \text{and}\ \ 
  (\theta, v) &\geq 0, \label{eq.kkt4}
  \end{align}
\esubequations
where $(u,v) \in \R{} \times \R{\ell}$ are Lagrange multipliers and $\Theta \equiv \operatorname{diag}(\theta)$.  For future use, we also introduce the notation $V \equiv \operatorname{diag}(v)$ for the nonnegativity multipliers.

The KKT conditions are difficult to solve directly due to the complementarity conditions represented by \eqref{eq.kkt3}--\eqref{eq.kkt4}.  Thus, as in a standard primal-dual interior-point method, rather than solve~\eqref{eq.kkt} directly, we introduce a centering term that, when employed, promotes the computation of steps toward or along the \emph{central path}.  This path, parameterized by the \emph{barrier parameter} $\mu \in (0,\infty)$, is defined by the set of points satisfying \eqref{eq.kkt1}, \eqref{eq.kkt2}, and \eqref{eq.kkt4} along with $\Theta v = \mu \ones$.

Given a solution estimate $(\theta_j,u_j,v_j)$, our method commences by computing a \emph{predictor} step corresponding to $\mu = 0$ 
as a solution of the linear system
\bequation\label{eq.lin_sys_pred}
  \bbmatrix
    -Q & A^T & I \\
    A & 0 & 0 \\
    V_j & 0 & \Theta_j
  \ebmatrix
  \bbmatrix
    \Delta \theta_j^{\text{pred}} \\
    \Delta u_j^{\text{pred}} \\
    \Delta v_j^{\text{pred}}
  \ebmatrix
  =
  \bbmatrix
    Q \theta_j + c - A^T u_j - v_j \\
    b - A \theta_j \\
    -\Theta_j v_j
  \ebmatrix
  =:
  \bbmatrix
    r_j^d \\
    r_j^p \\
    r_j^c
  \ebmatrix.
\eequation
The residual vectors $r^p$ and $r^d$ correspond to \textit{primal} and \textit{dual} infeasibility, respectively, whereas the residual vector $r^c$ corresponds to the \textit{complementarity} condition.  The linear system \eqref{eq.lin_sys_pred} can be solved directly, or it can be solved by first solving
\bequation\label{eq.lin_sys_pred_reduced}
  \(
  \underbrace{
  \bbmatrix
    -Q & A^T \\
    A & 0
  \ebmatrix
  }_{\text{constant w.r.t.~$j$}}
  +
  \bbmatrix
    - \Theta_j^{-1} V_j & 0 \\
    0 & 0
  \ebmatrix
  \)
  \bbmatrix
    \Delta \theta_j^{\text{pred}} \\
    \Delta u_j^{\text{pred}} \\
  \ebmatrix
  =
  \bbmatrix
    r_j^d - \Theta_j^{-1} r_j^c \\
    r_j^p
  \ebmatrix,
\eequation
then setting the remaining solution component as $\Delta v_j^{\text{pred}} \gets \Theta_j^{-1}(r_j^c - V_j \Delta \theta_j^{\text{pred}})$.  Once the predictor step $(\Delta \theta_j^{\text{pred}}, \Delta u_j^{\text{pred}}, \Delta v_j^{\text{pred}})$ has been computed, our method next computes step sizes $(\alpha_{\theta,j}^{\text{pred}},\alpha_{u,j}^{\text{pred}},\alpha_{v,j}^{\text{pred}})$ following a strategy based on that in~\cite{CurtNoce07}, which is discussed further in the next paragraph.  Then, given a current barrier parameter value $\mu_j \gets \theta_j^Tv_j/\ell$ along with the predictor step and its corresponding step sizes, the method next computes, as in a standard Mehrotra predictor-corrector strategy~\cite{Mehr1992}, the scaling factor $\zeta_j \in \R{}_{>0}$ by
\bequation\label{eq.scaling}
  \zeta_j \gets \(\frac{1}{\mu_j} \cdot \frac{(\theta_j + \alpha_{\theta,j}^{\text{pred}} \Delta \theta_j^{\text{pred}})^T(v_j + \alpha_{v,j}^{\text{pred}} \Delta v_j^{\text{pred}})}{\ell}\)^3.
\eequation
(In our implementation, we safeguard this choice by setting a larger value for $\zeta_j$, if necessary, to ensure that $\zeta_j\mu_j \geq 10^{-12}$.)  The method next computes a \emph{corrector} step based on this scaling factor.  Equivalently, the method computes the full search direction (i.e., the predictor step plus the corrector step) by solving
\bequation\label{eq.lin_sys}
  \bbmatrix
    -Q & A^T & I \\
    A & 0 & 0 \\
    V_j & 0 & \Theta_j
  \ebmatrix
  \bbmatrix
    \Delta \theta_j \\
    \Delta u_j \\
    \Delta v_j
  \ebmatrix
  =
  \bbmatrix
    Q \theta_j + c - A^T u_j - v_j \\
    b - A \theta_j \\
    -\Theta_j v_j + \zeta_j \mu_j \ones
  \ebmatrix
  =:
  \bbmatrix
    r_j^d \\
    r_j^p \\
    r_j^c + \zeta_j \mu_j \ones
  \ebmatrix.
\eequation
Like for the predictor step, the linear system \eqref{eq.lin_sys} can be solved directly, or it can be solved by employing a reduced linear system similar to~\eqref{eq.lin_sys_pred_reduced}.  In any case, our implementation exploits the fact that the predictor and corrector steps are computed through linear systems involving the same matrix, which allows our implementation to re-use the LDL factorization from the computation of the predictor step in order to compute the corrector step.  Once the full step has been computed, the method computes corresponding step sizes $(\alpha_{\theta,j}, \alpha_{u,j}, \alpha_{v,j})$, again using the method proposed in~\cite{CurtNoce07} (see the next paragraph).  Finally, the solution estimate is updated and the method proceeds to the next iteration.

The search-direction-computation strategy described in the previous paragraph has been employed to great effect for solving linear and quadratic optimization problems with interior-point methods.  That said, a critical practical aspect of our proposed method is the manner in which step sizes are chosen.  For this purpose, we have found impressive performance to be offered by the strategy proposed in \cite{CurtNoce07}, tailored for our setting.  This approach starts by choosing maximum step sizes such that a \emph{fraction-to-the-boundary} rule is enforced; see line~\ref{step.ftb} in Algorithm~\ref{alg.step_sizes}.  Afterwards, separate primal and dual step sizes are computed in order to optimize the overall reduction in the merit function $\phi : \R{\ell} \times \R{} \times \R{\ell}$ defined by
\bequationNN
  \phi(\theta, u, v) = \|A\theta-b\|_2^2 + \|A^Tu + v - Q\theta - c\|_2^2 + \theta^Tv.
\eequationNN
This is done by searching along a two-segment path from $(0,0)$ to the maximum step sizes allowed by the fraction-to-the-boundary rule.  We refer the reader to \cite{CurtNoce07} for further details.  For our purposes here, we simply state that the desired step sizes can be computed by solving two 2-dimensional convex quadratic optimization problems, the necessary computations for which are negligible.

A complete description of our entire interior-point method is stated as Algorithm~\ref{alg.IPM} along with the step-size-selection subroutine stated in Algorithm~\ref{alg.step_sizes}.  Besides tailoring the aforementioned strategies and linear algebra subroutines to the specific formulation in \eqref{eq.dual_subproblem}, the initial point provided to Algorithm~\ref{alg.IPM} has been tuned carefully to obtain improved practical performance.  We discuss this initialization further in the next section along with other implementation details of NonOpt.

\balgorithm
  \caption{Interior Point Algorithm for Solving Problem \eqref{eq.quad_prob}}
  \label{alg.IPM}
  \balgorithmic[1]
    \State choose initial point $(\theta_1, u_1, v_1) \in \R{\ell}_{>0} \times \R{} \times \R{\ell}_{>0}$ and termination tolerance $\epsilon \in \R{}_{>0}$
    \For{\textbf{all} $j \in \N{}$}
      \State set $(r_j^p,r_j^d,r_j^c)$ as defined in \eqref{eq.lin_sys_pred}
      \If{$\max\{\|r_j^p\|_\infty, \|r_j^d\|_\infty, \|r_j^c\|_\infty\} \leq \epsilon$}
        \State \textbf{terminate} and \Return $(\theta_j, u_j, v_j)$
      \EndIf
      \State compute $(\Delta \theta_j^{\text{pred}}, \Delta u_j^{\text{pred}}, \Delta v_j^{\text{pred}})$ by solving the linear system in \eqref{eq.lin_sys_pred}
      \State compute $(\alpha_{\theta,j}^{\text{pred}}, \alpha_{u,j}^{\text{pred}}, \alpha_{v,j}^{\text{pred}})$ by Algorithm~\ref{alg.step_sizes} for $(\Delta \theta_j^{\text{pred}}, \Delta u_j^{\text{pred}}, \Delta v_j^{\text{pred}})$
      \State set $\mu_j \gets \theta_j^Tv_j/\ell$
      \State set $\zeta_j$ by \eqref{eq.scaling}
      \State compute $(\Delta \theta_j, \Delta u_j, \Delta v_j)$ by solving the linear system in \eqref{eq.lin_sys}
      \State compute $(\alpha_{\theta,j}, \alpha_{u,j}, \alpha_{v,j})$ by Algorithm~\ref{alg.step_sizes} for $(\Delta \theta_j, \Delta u_j, \Delta v_j)$
      \State set $\theta_{j+1} \gets \theta_j + \alpha_{\theta,j} \Delta \theta_j$, $u_{j+1} \gets u_j + \alpha_{u,j} \Delta u_j$, and $v_{j+1} \gets v_j + \alpha_{v,j} \Delta v_j$
    \EndFor
  \ealgorithmic
\ealgorithm

\balgorithm
  \caption{Step-size computation subroutine for Algorithm~\ref{alg.IPM}}
  \label{alg.step_sizes}
  \balgorithmic[1]
    \Require search direction, call it $(\Delta \theta_j, \Delta u_j, \Delta v_j)$
    \State get current iterate $(\theta_j, u_j, v_j)$ and residuals $(r_j^p,r_j^d)$ from Algorithm~\ref{alg.IPM}
    \State set $\overline\alpha_{\theta,j} \gets \left[ \max_{i \in [\ell]} \left\{ 1, -\frac{\Delta \theta_j^{(i)}}{\beta \theta_j^{(i)}} \right\} \right]^{-1}$ and $\overline\alpha_{v,j} \gets  \left[ \max_{i \in [\ell]} \left\{ 1, -\frac{\Delta v_j^{(i)}}{\beta v_j^{(i)}} \right\} \right]^{-1}$ \label{step.ftb}
    \State set $\overline\alpha_j \gets \min \{ \overline\alpha_{\theta,j}, \overline\alpha_{v,j} \}$
    \State set
    \bequationNN
      r_j \gets \bbmatrix r_j^p \\ r_j^d \ebmatrix,\ s_j \gets \bbmatrix -A \\ Q \ebmatrix \Delta \theta_j,\ o_j \gets \bbmatrix 0 \\ -A^T \ebmatrix \Delta u_j,\ \text{and}\ p_j \gets \bbmatrix 0 \\ -I \ebmatrix \Delta v_j
    \eequationNN
    \State compute $(\alpha_{\theta,j,1},\alpha_{u,j,1})$ by solving the two-dimensional convex QP given by
    \begin{align*}
      \min_{(\alpha_\theta, \alpha_u) \in \R{} \times \R{}} &\ \half \bbmatrix \alpha_\theta \\ \alpha_u \ebmatrix \bbmatrix (s_j + p_j)^T (s_j + p_j) + \Delta \theta_j^T \Delta v_j & (s_j + p_j)^T o_j \\ (s_j + p_j)^T o_j & o_j^T o_j \ebmatrix \bbmatrix \alpha_\theta \\ \alpha_u \ebmatrix \\
      &\ + \bbmatrix r_j^T (s_j + p_j) + \thalf (\Delta \theta_j^Tv_j + \theta_j^T \Delta v_j) \\ r_j^T o_j \ebmatrix^T \bbmatrix \alpha_\theta \\ \alpha_u \ebmatrix \\
      \st &\ 0 \leq \alpha_\theta \leq \overline\alpha_j
    \end{align*}
    \State set $\alpha_{v,j,1} \gets \alpha_{\theta,j,1}$
    \If{$\overline\alpha_{\theta,j} \leq \overline\alpha_{v,j}$}
      \State set $\alpha_{\theta,j,2} \gets \overline\alpha_{\theta,j}$
      \State compute $(\alpha_{u,j,2},\alpha_{v,j,2})$ by solving the two-dimensional convex QP given by
      \begin{align*}
        \min_{(\alpha_u, \alpha_v) \in \R{} \times \R{}} &\ \half \bbmatrix \alpha_u \\ \alpha_v \ebmatrix \bbmatrix o_j^T o_j & o_j^T p_j \\ o_j^T p_j & p_j^T p_j \ebmatrix \bbmatrix \alpha_u \\ \alpha_v \ebmatrix \\
        &\ + \bbmatrix r_j^T o_j + \overline\alpha_{\theta,j} s_j^T o_j \\ r_j^T p_j + \thalf \theta_j^T \Delta v_j + \overline\alpha_{\theta,j} (s_j^T p_j + \thalf \Delta \theta_j^T \Delta v_j) \ebmatrix^T \bbmatrix \alpha_u \\ \alpha_v \ebmatrix \\
        \st &\ \overline\alpha_j \leq \alpha_v \leq \overline\alpha_{v,j}
      \end{align*}
    \Else
      \State compute $(\alpha_{\theta,j,2},\alpha_{u,j,2})$ by solving the two-dimensional convex QP given by
      \begin{align*}
        \min_{(\alpha_\theta, \alpha_u) \in \R{} \times \R{}} &\ \half \bbmatrix \alpha_\theta \\ \alpha_u \ebmatrix \bbmatrix s_j^T s_j & s_j^T o_j \\
s_j^T o_j & o_j^T o_j \ebmatrix \bbmatrix \alpha_\theta \\ \alpha_u \ebmatrix \\
        &\ + \bbmatrix r_j^T s_j + \thalf \Delta \theta_j^T v_j + \overline\alpha_{v,j} (s_j^T p_j + \thalf \Delta \theta_j^T \Delta v_j) \\ r_j^T o_j + \overline\alpha_{v,j} o_j^T p_j \ebmatrix^T \bbmatrix \alpha_\theta \\ \alpha_u \ebmatrix \\
        \st &\ \overline\alpha_j \leq \alpha_\theta \leq \overline\alpha_{\theta,j}
      \end{align*}
      \State set $\alpha_{v,j,2} \gets \overline\alpha_{v,j}$
    \EndIf
    \State set $\chi \gets \displaystyle \arg\min_{\chi \in \{1,2\}} \phi(\theta_j + \alpha_{\theta,j,\chi} \Delta \theta_j, u_j + \alpha_{u,j,\chi} \Delta u_j, v_j + \alpha_{v,j,\chi} \Delta v_j)$
    \State set $\alpha_{\theta,j} \gets \alpha_{\theta,j,\chi}$, $\alpha_{u,j} \gets \alpha_{u,j,\chi}$, and $\alpha_{v,j} \gets \alpha_{v,j,\chi}$
    \State \Return $(\alpha_{\theta,j}, \alpha_{u,j}, \alpha_{v,j})$
  \ealgorithmic
\ealgorithm

\section{Implementation Details}\label{sec.implementation}

The performance of a gradient-based algorithm for solving \emph{smooth} nonconvex continuous optimization problems is sensitive to input parameter choices in general.  The performance of an algorithm for solving \emph{nonsmooth} nonconvex continuous optimization problems is arguably even more sensitive.  We have certainly found this to be the case for NonOpt.  Fortunately, based on its default input parameter settings, most of the behaviors of the implemented algorithms in NonOpt seem to be similar to that of an unadulterated quasi-Newton scheme with a weak-Wolfe line search, which is generally quite effective even when solving locally Lipschitz problems that involve nonconvexity and/or nonsmoothness \cite{LewiOver2013}.  This can be seen in our experiments in the next section, where one finds (see \S\ref{sec.speed}--\S\ref{sec.accuracy} in particular) that the gradient-based direction computation strategy (indicated by \texttt{G}) generally performs quite well.  (However, it is also clear from our experiments in the next section that, for certain problems, the gradient-based strategy is inferior to the gradient-sampling and/or proximal-bundle strategy.)  NonOpt employs the default input parameters of $c \gets 10^{-10}$, $\eta \gets 10^{-8}$, and $\psi \gets 10^8$ so that, for the most part, an employed algorithm follows such an unadulterated line-search-BFGS scheme.  That being said, if a user wants to push an algorithm to obtain a higher-accuracy solution when nonsmoothness is present, then one might consider alternative values of these inputs (e.g., setting $(\eta,\psi) \gets (10^{-2},10^2)$) to force each local model of the objective (recall~\eqref{eq.primal_subproblem}) to rely more heavily on the information in the point set $\{(x_{k,j},f_{k,j},g_{k,j})\}_{j\in[m]}$ than the information embedded in the Hessian approximation $H_k$.  The user should understand, however, that such a choice often leads to the algorithm requiring a much higher computational cost per iteration.  This, at least, is true based on our observations that such a choice leads to shorter steps, more accumulation of points in the point sets, and thus higher subproblem costs.

Any user of NonOpt should also be made aware that the strategy for choosing~$\epsilon_1 \in \R{}_{>0}$ and setting the remaining elements of the generated sequence $\{\epsilon_k\}$ play an absolutely critical role in the behavior of the software.  To describe this strategy, let us being by saying that NonOpt commences any solve by scaling the objective function in order to ensure that $\|g_1\|_\infty$ is at most $10^2$.  Letting $\ftilde$ denote the pre-scaled objective function and $\gtilde_1 \in \gengradient f(x_1)$, this is accomplished by setting the objective function to be considered in the remainder of the solve as
\bequationNN
  f \gets \texttt{scale} \cdot \ftilde,\ \ \text{where}\ \ \texttt{scale} \gets \min\{1, 10^2 \|\gtilde_1\|_\infty^{-1} \}.
\eequationNN
After this scaling, by default, the software sets $\epsilon_1 \gets \max\{10^{-2}, 10^{-1} \|g_1\|_\infty\}$.  As for setting $\epsilon_{k+1} \in (0,\epsilon_k]$ for all $k \in \N{}$, the ideal case is for the software to set $\epsilon_{k+1} < \epsilon_k$ if and only if (see Theorems 3.3--3.4 in \cite{CurtRobiZhou20})
\bequationNN
  \max\{\|d_k\|_\infty, \|G_k\omega_k\|_\infty, \|G_k\omega_k + \gamma_k\|_\infty\} \leq \epsilon_k.
\eequationNN
If and when this condition is satisfied, the software sets $\epsilon_{k+1} \gets 10^{-1} \epsilon_k$.  However, for practical purposes, the software will also set $\epsilon_{k+1} \gets 10^{-1} \epsilon_k$ whenever recent improvements in the objective function are small.  We discuss this further in \S\ref{sec.speed}.

The trust region in the primal subproblem~\eqref{eq.primal_subproblem} can be an effective tool for obtaining a high-accuracy solution.  Indeed, setting $\{\delta_k\}$ to relatively small values has a similar effect as choosing $\eta$ and $\psi$ each closer to 1, as discussed earlier in this section.  However, a disadvantage of setting $\delta_k$ too small is that it causes the subproblem solve to be more expensive in general, which can be especially detrimental when solving large-scale problems.  Therefore, by default, the software sets $\delta_1 \gets \max\{10^{-1}, 10^{10} \|g_1\|_\infty\}$, then, for each $k \in \N{}$, sets $\delta_{k+1} \gets \delta_k$ unless $\epsilon_{k+1} < \epsilon_k$, in which case it sets $\delta_{k+1} \gets 10^{-1} \delta_k$.  These default settings make it unlikely for the trust-region constraint ever to be active.  Consequently, for the implemented interior-point solver, the software first attempts to solve \eqref{eq.quad_prob} without introducing the variables $(\sigma,\rho)$.  If the resulting solution $\tilde\omega$ yields $\|W_kG_k\tilde\omega\|_\infty \leq \delta_k$, then the software sets $(\omega_k,\sigma_k,\rho_k) \gets (\tilde\omega,0,0)$.  Otherwise, it re-solves \eqref{eq.quad_prob} with all variables in order to compute $(\omega_k,\sigma_k,\rho_k)$.  If $n$ is very large, then this latter occurrence can unfortunately lead to a very expensive subproblem solve.

As previously mentioned, the subproblem solves are the most computationally expensive aspect of each iteration of the algorithms in NonOpt, except possibly for the function and generalized gradient evaluations.  This means, in turn, both the initialization strategies for the subproblem solvers and the linear algebra subroutines employed in the solvers are critically important.  The dual active-set subproblem solver for \eqref{eq.dual_subproblem} initializes its solution estimate to have $\omega_j = 1$ for some $j \in [m]$ and $\omega_j = 0$ for all other elements of $[m]$, as well as $\gamma = 0$.  The index $j$ is chosen as the index yielding the smallest value of $g_{k,j}^TW_kg_{k,j}$ over $j \in [m]$.  The interior-point solver is initialized with $\omega_j \gets 1/m$ for all $j \in [m]$, so that it is also initialized with a feasible solution for \eqref{eq.dual_subproblem}.  Moreover, with respect to $\mu_0 \gets 5 \times 10^{-1}$ (a value chosen based on extensive experimentation), it initializes the components of $v$ corresponding to $\omega$, say, represented by a vector $v^\omega$, as $\mu / \omega_j$ for all $j \in [m]$.  This ensures that the initial point is complementary in the sense that $\omega_jv_j^\omega = \mu$ for all $j \in [m]$.  As for $(\sigma,\rho)$, if they are introduced, their components are initialized based on the values of $W_kG_k\omega$ for the initialized $\omega$.  If $[W_kG_k\omega]_i < -\delta_k$, then
\bequationNN
  \sigma_i \gets \max\{10^{-1}, \delta_k - [W_kG_k\omega]_i\},\ \rho_i \gets 10^{-1}, v_i^\sigma \gets \mu_0 / \sigma_i,\ \text{and}\ v_i^\rho \gets \mu_0 / \rho_i,
\eequationNN
whereas if $[W_kG_k\omega]_i > \delta_k$, then
\bequationNN
  \sigma_i \gets 10^{-1},\ \rho_i \gets \max\{10^{-1}, -\delta_k - [W_kG_k\omega]_i\}, v_i^\sigma \gets \mu_0 / \sigma_i,\ \text{and}\ v_i^\rho \gets \mu_0 / \rho_i.
\eequationNN
The use of the value $10^{-1}$ (chosen based on extensive experimentation) ensures that the initial point is an interior point.  Careful selection of the initial point led to a significant reduction in CPU time for our experiments in the next section.  We regularly saw that the number of iterations was around 25\% of those required when the solve was conducted without a carefully selected initial point.

The interior-point subproblem solver in NonOpt relies on the LAPACK subroutines \texttt{dsysv} (to factor the matrix for computing the predictor step) and \texttt{dsytrs} (for re-using the factorization to compute the full interior-point step).  The dual active-set subproblem solver generally computes and updates its own Cholesky factorizations based on iterative additions and deletions; see \cite{Kiwi1986}.  In such cases, linear systems are solved using LAPACK's \texttt{dtrsv} subroutine.  However, if an update of a factorization through an addition is suspected to lead to numerical error, then a Choleksy factorization of a matrix is computed from scratch using LAPACK's \texttt{dpotf2} subroutine.  An ambitious user could replace these subroutines with ones, say, from the Harwell Subroutine Library\footnote{\href{https://www.hsl.rl.ac.uk/}{\texttt{https://www.hsl.rl.ac.uk/}}} or the Pardiso Solver Project.\footnote{\href{https://panua.ch/pardiso/}{\texttt{https://panua.ch/pardiso/}}}

\section{Numerical Experiments}\label{sec.numerical}

We ran four sets of experiments to demonstrate various features of the NonOpt software as well as its overall effectiveness in solving challenging locally Lipschitz minimization problems.  Previous articles about the developments of algorithmic strategies that have led to the NonOpt software considered comparisons with other software packages, such as LMBM \cite{Karm2007}.  However, since NonOpt performed favorably in these comparisons, LMBM is no longer under active development, and we are not aware of any other open-source software packages under active development for general-purpose locally Lipschitz minimization, our experiments in this paper focus exclusively on experiments with NonOpt only.

Our first set of experiments involves comparing the performance of the dual active-set QP subproblem solver versus the newly implemented interior-point QP subproblem solver (described in detail in \S\ref{sec.subproblem}).  These experiments show that the interior-point QP solver is a very valuable addition to the software.

Our next two sets of experiments aim to demonstrate the versatility of the NonOpt software when solving a diverse set of locally Lipschitz test problems.  The two sets of experiments can be motivated as follows.  First, in certain settings, a user may desire for the software to be run for speed, say, to obtain a decent solution relatively quickly.  Therefore, we have one set of experiments where the user-defined parameters lead to relatively fast termination.  On the other hand, in other settings, a user may desire for the software to be run for accuracy, say, to obtain a higher-quality solution even if the computational time required is significantly higher.  Therefore, we have another set of experiments (on the same set of test problems as the previous set of experiments) where the user-defined parameters lead to longer solution times, but better final solutions.

Finally, we present the results of experiments with a few large-scale problems from the field of image processing.  The purpose of these experiments is to show that NonOpt is capable of solving problems with hundreds-of-thousands of variables in a computationally efficient manner, and that it can transition seamlessly from solving convex to solving nonconvex locally Lipschitz minimization problems.

All of the experimental results that are provided in this section were obtained by running NonOpt version~2.0 on a Macbook Pro running Sequoia 15.3.1 with an Apple M1 Pro chip and 16GB of memory.

\subsection{Experiments with the QP subproblem solvers}

The purpose of our first set of experiments is to demonstrate the relative performance of the newly implemented interior-point QP subproblem solver as compared to the previously implemented dual active-set QP solver.  Let us refer to these as the IPM and DAS solvers, respectively.  The DAS solver is computationally efficient when solving relatively small-scale subproblems and in certain instances when a subproblem is solved that is closely related to the most recently solved subproblem, such as when the software has solved an instance of \eqref{eq.dual_subproblem}, then requests for a solution of a related instance for which the only differences are that the matrix $G_k$ has been augmented with additional columns and, correspondingly, the vector $b_k$ has been augmented with additional components.  However, especially when solving large-scale subproblems, the IPM solver can be much faster, as is typical when one compares interior-point methods with active-set-type methods.

For these experiments, we randomly generated a large number of QP subproblems of increasing size.  The generation of each QP proceeded as explained in the bullets below, which essentially amounts to generating values to satisfy the KKT conditions in~\eqref{eq.kkt} in order to correspond with a prescribed choice of the optimal solution of the primal subproblem~\eqref{eq.primal_subproblem}.  In the following description, we refer to the quantities that one needs to define for an instance of subproblem~\eqref{eq.dual_subproblem}, but drop the iteration $k$ index.  The generation procedure for each instance is given problem size parameters $n \in \N{}$ and $m \in \N{}$.  The values for these parameters that we considered are discussed after the following description of the generation procedure.

\bitemize
  \item First, $W \in \R{n \times n}$ was set to the identity matrix.  This choice is reasonable since, as far as the subproblem solvers are concerned, the positive-definite matrix $W$ essentially only transforms the columns of $G$.  We expect the relative performance of the subproblem solvers to be the same for other choices of $W$.
  \item Second, $\delta \in \R{}_{>0}$ was set to 1 and an option was set to decide whether the optimal solution of the primal QP, call it $d_*$, would have (a) $d_* = 0$, (b)~$[d_*]_i = \delta$ for all $i \in [n/2]$ and $[d_*] = 0$ for all $i \in \{n/2+1,\dots,n\}$, or (c)~$d_* = \delta \ones$.  (We generated an equal number of each of these types of QPs for each problem size.)  Corresponding to these cases, an unconstrained minimizer of the objective function of~\eqref{eq.primal_subproblem}, call it $d_{unc}$, was set to have (a) $d_{unc} = 0$, (b) $[d_{unc}]_i = 2\delta$ for all $i \in [n/2]$ and $[d_*] = 0$ for all $i \in \{n/2+1,\dots,n\}$, or (c) $d_{unc} = 2 \delta \ones$.
  \item Third, $G \in \R{n \times m}$ and $\omega_* \in \R{m}$ were generated to satisfy $\ones^T\omega_* = 1$, $\omega_* \geq 0$, and $G\omega_* = -d_{unc}$.  This was done by first randomly generating $\Ghat \in \R{n \times (n-1)}$ with each component drawn from a standard normal distribution and randomly generating $\hat\omega \in \R{n-1}$ with each component drawn from a uniform distribution over $[0,1]$.  Then, $\Ghat$ was padded with the additional column $-d_{unc} - \Ghat\hat\omega$ and~$\hat\omega$ was padded with an additional component equal to 1.  The next step was to compute $\zeta \gets \|\hat\omega\|_1$, normalize $\hat\omega \gets \hat\omega/\zeta$ so that $\ones^T\hat\omega = 1$, and correspondingly set $\Ghat \gets \zeta \Ghat$.  Finally, $\Ghat$ was padded with an additional $m-n$ columns with components drawn from a standard normal distribution and~$\hat\omega$ was padded with additional 0s.  This procedure guarantees $\ones^T\omega_* = 1$, $\omega_* \geq 0$, and $G\omega_* = -d_{unc}$.
  \item Fourth, $q \gets -(d_* + WG\omega_*)$ was computed in order to set $\rho_* \gets \max\{-q,0\}$ and $\sigma_* \gets \max\{q,0\}$, where the $\max$ is taken component-wise in each case.  (This guarantees that $d_* = -W(G\omega_* + \sigma_* - \rho_*)$, where $\sigma_* \geq 0$ and $\rho_* \geq 0$.)
  \item Fifth, $u$ was set to an arbitrary number.  (For our experiments in this section, without loss of generality we simply set $u \gets 5$.)
  \item Sixth, all that remains is to set $b \in \R{m}$.  Toward this end, a vector $v_\omega \in \R{m}$ was generated to have $[v_\omega]_i = 0$ for all $i \in [m]$ such that $[\omega_*]_i > 0$, and all other components drawn from a uniform distribution over [0,1].  Finally, $b \gets G^TW(G\omega_* + \sigma_* - \rho_*) + u\ones - v_\omega$.
\eitemize

One can verify that, with the procedure above, the generated QP is such that the corresponding unique optimal solution of \eqref{eq.primal_subproblem} is the prescribed vector~$d_*$.  Overall, we generated 900 QP test problems.  We generated equal numbers of problems for the values of $n \in \{100,200,300,400,5000,600,700,800,900,1000\}$, and for each~$n$ we generated equal numbers of problems for the values of $m \in \{n+1,1.5n,2n\}$.  As previously mentioned, for each $(n,m)$ we generated equal numbers of problems for the three prescribed choices of $d_*$, and for each setting of $(n,m,d_*)$ we generated 10 different problems, each using a different seed for the random number generator.  This led to the 900 different QP subproblems.

\begin{figure}[ht]
  \centering
  \includegraphics[width=0.75\textwidth,clip=true,trim=90 00 90 10]{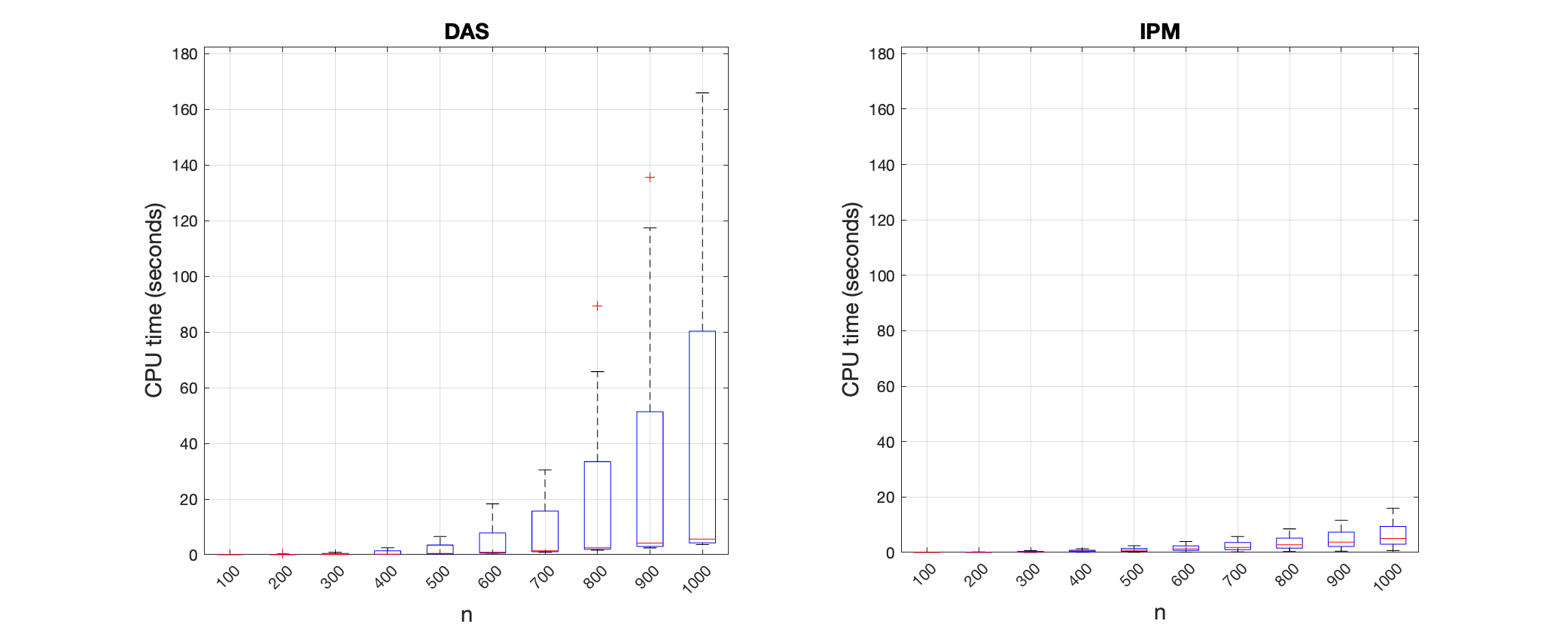}
  \caption{DAS vs. IPM, all problems}
  \label{fig.das_ipm_0}
\end{figure}

CPU time results for all runs combined are provided in Figure~\ref{fig.das_ipm_0}.  As $n$ increases, it is clear that IPM scales much better than DAS in terms of solution time.  That said, the solution time for DAS is highly variable, so it is appropriate to parse the results further to understand if there are scenarios or settings in which DAS may be more competitive for IPM, even when $n$ increases.

\begin{figure}[ht]
  \centering
  \includegraphics[width=0.75\textwidth,clip=true,trim=90 00 90 10]{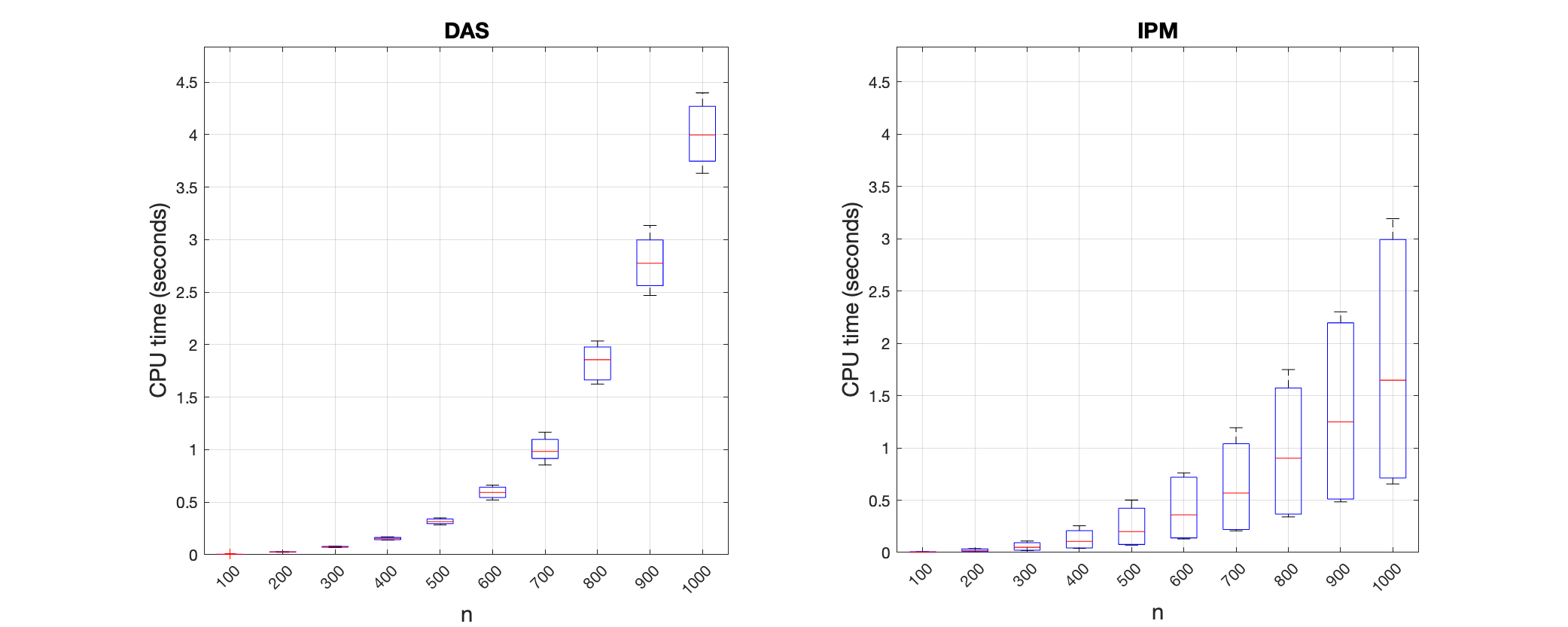}
  \caption{DAS vs. IPM, $d_* = 0$}
  \label{fig.das_ipm_1}
\end{figure}

\begin{figure}[ht]
  \centering
  \includegraphics[width=0.75\textwidth,clip=true,trim=90 00 90 10]{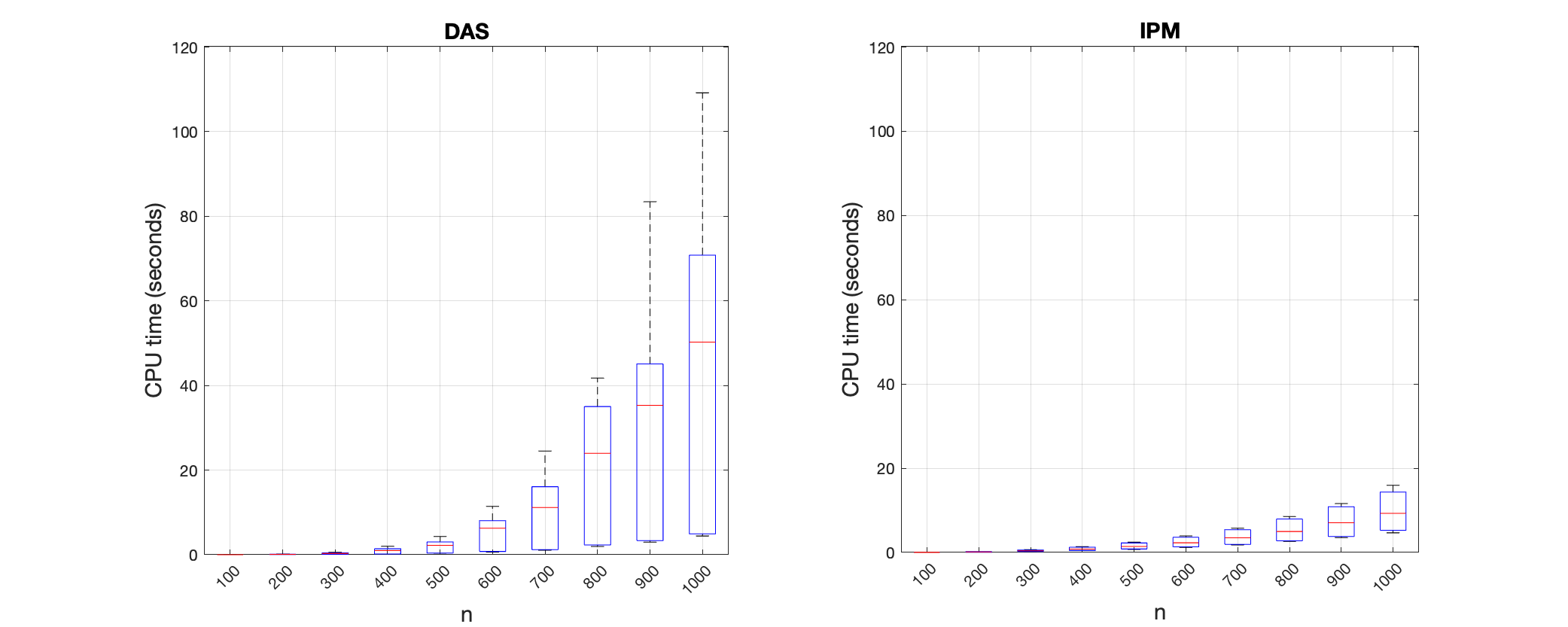}
  \caption{DAS vs. IPM, $[d_*]_i = \delta$ for all $i \in [n/2]$, 0 otherwise}
  \label{fig.das_ipm_2}
\end{figure}

\begin{figure}[ht]
  \centering
  \includegraphics[width=0.75\textwidth,clip=true,trim=90 00 90 10]{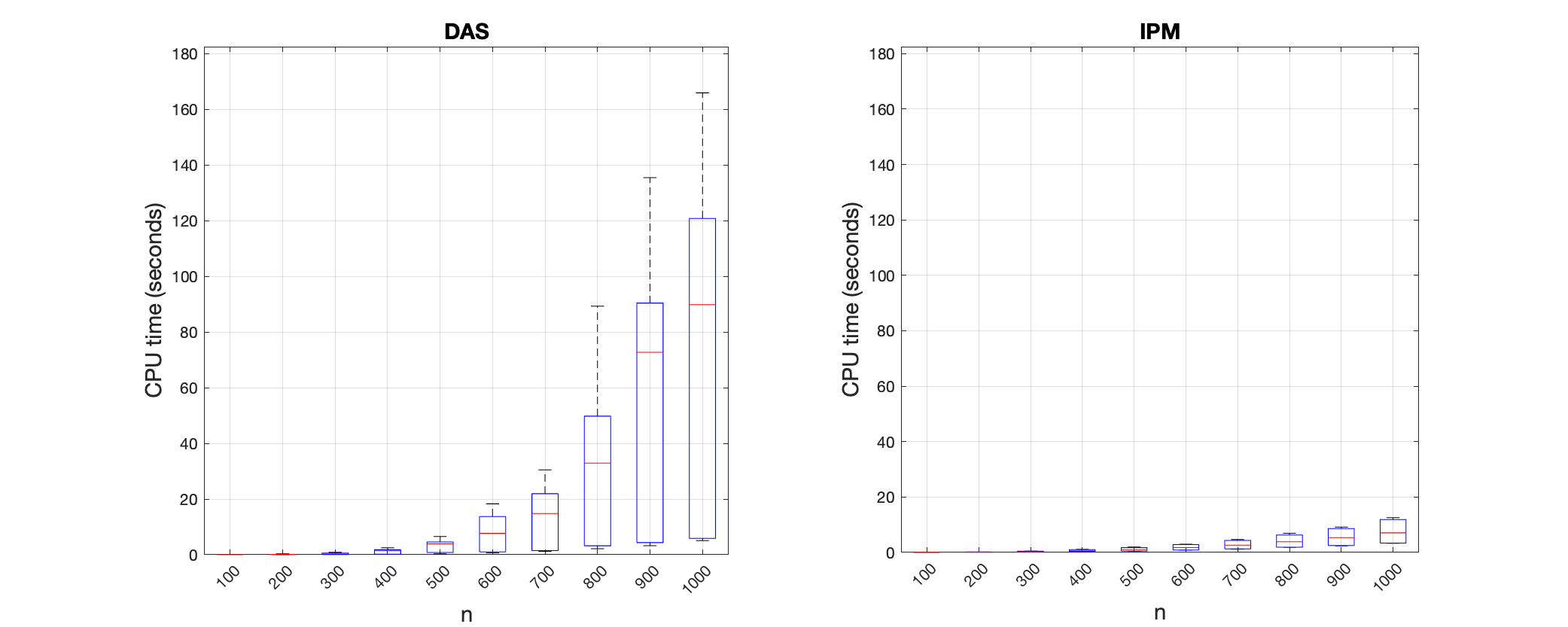}
  \caption{DAS vs. IPM, $d_* = \delta \ones$}
  \label{fig.das_ipm_3}
\end{figure}

Figures~\ref{fig.das_ipm_1}--\ref{fig.das_ipm_3} parse the results based on different $d_*$ values.  Through these results, one finds that DAS is more competitive with IPM when $d_* = 0$, which more generally suggests that DAS can be more competitive with IPM when the trust-region constraint is inactive.  However, even in this setting, IPM scales better than DAS as $n$ increases.  When the trust-region constraint is active, the IPM solution times are typically much faster than the DAS solution times.

\begin{figure}[ht]
  \centering
  \includegraphics[width=0.75\textwidth,clip=true,trim=90 00 90 10]{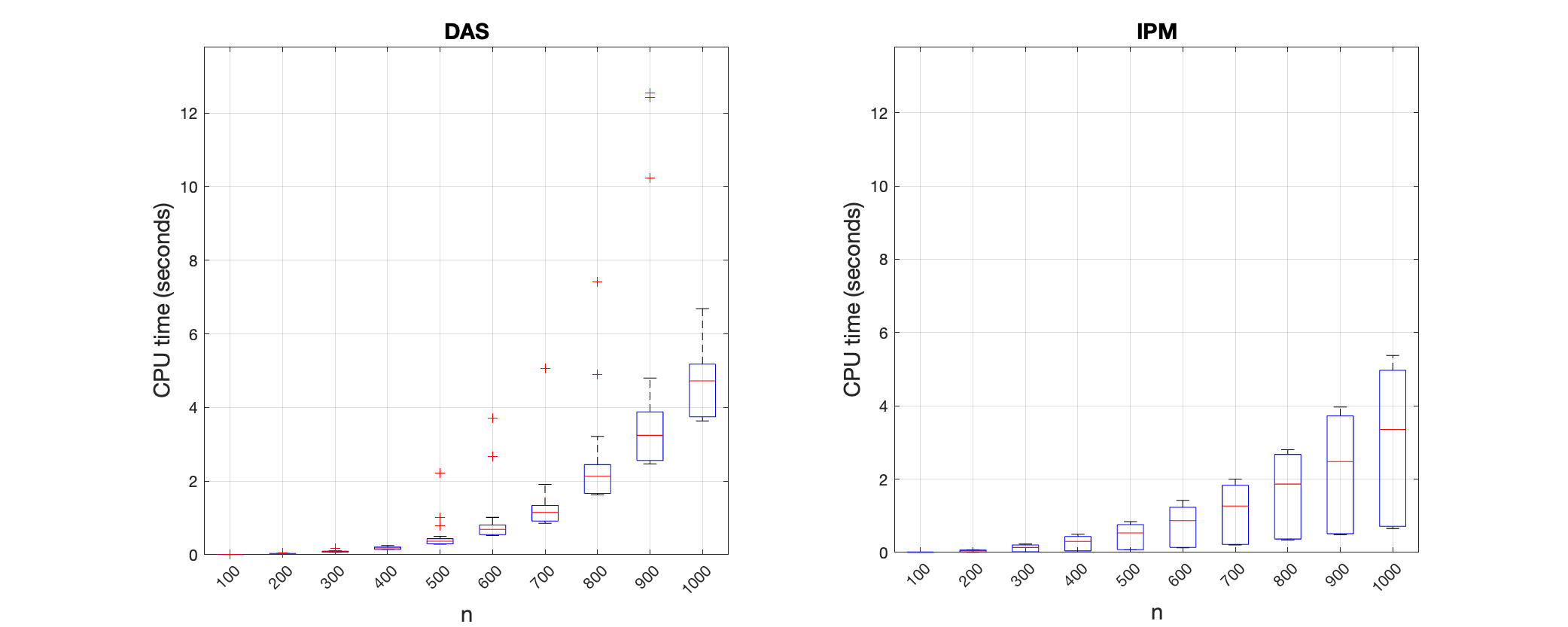}
  \caption{DAS vs. IPM, $m = n+1$}
  \label{fig.das_ipm_4}
\end{figure}

\begin{figure}[ht]
  \centering
  \includegraphics[width=0.75\textwidth,clip=true,trim=90 00 90 10]{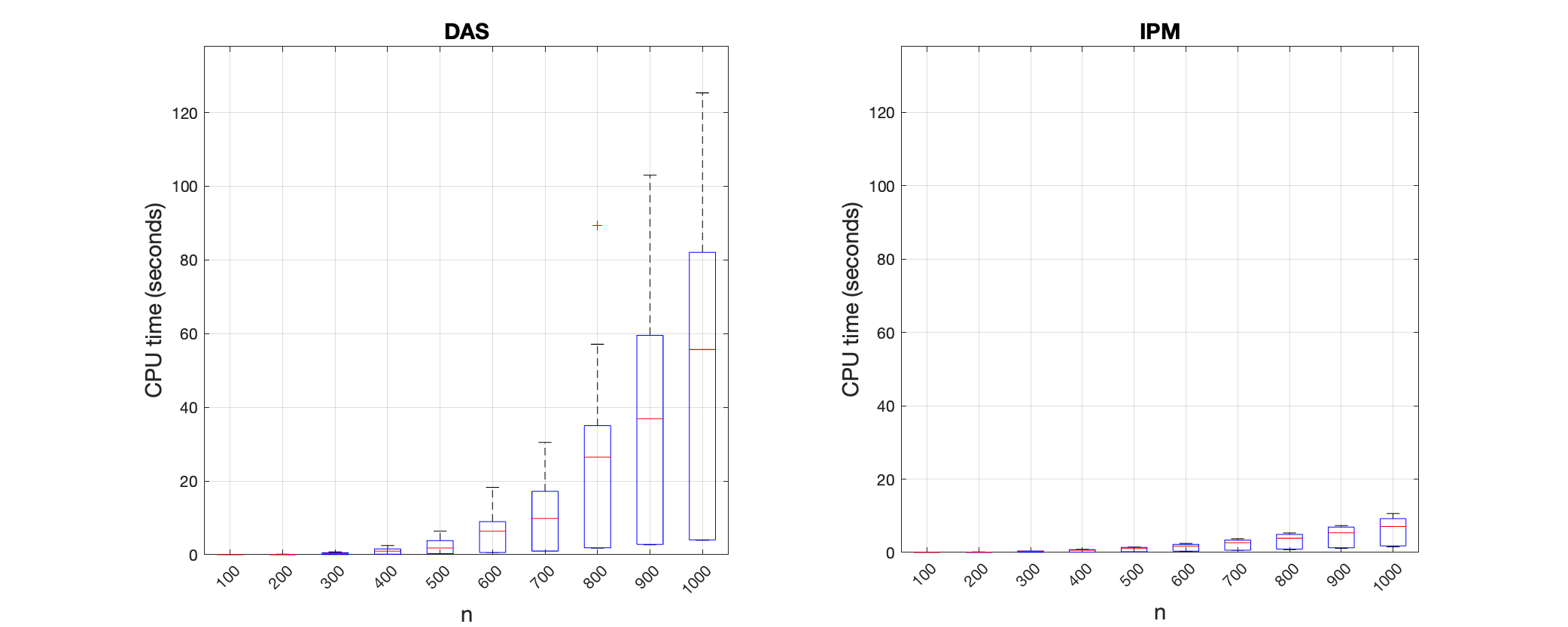}
  \caption{DAS vs. IPM, $m = 1.5n$}
  \label{fig.das_ipm_5}
\end{figure}

\begin{figure}[ht]
  \centering
  \includegraphics[width=0.75\textwidth,clip=true,trim=90 00 90 10]{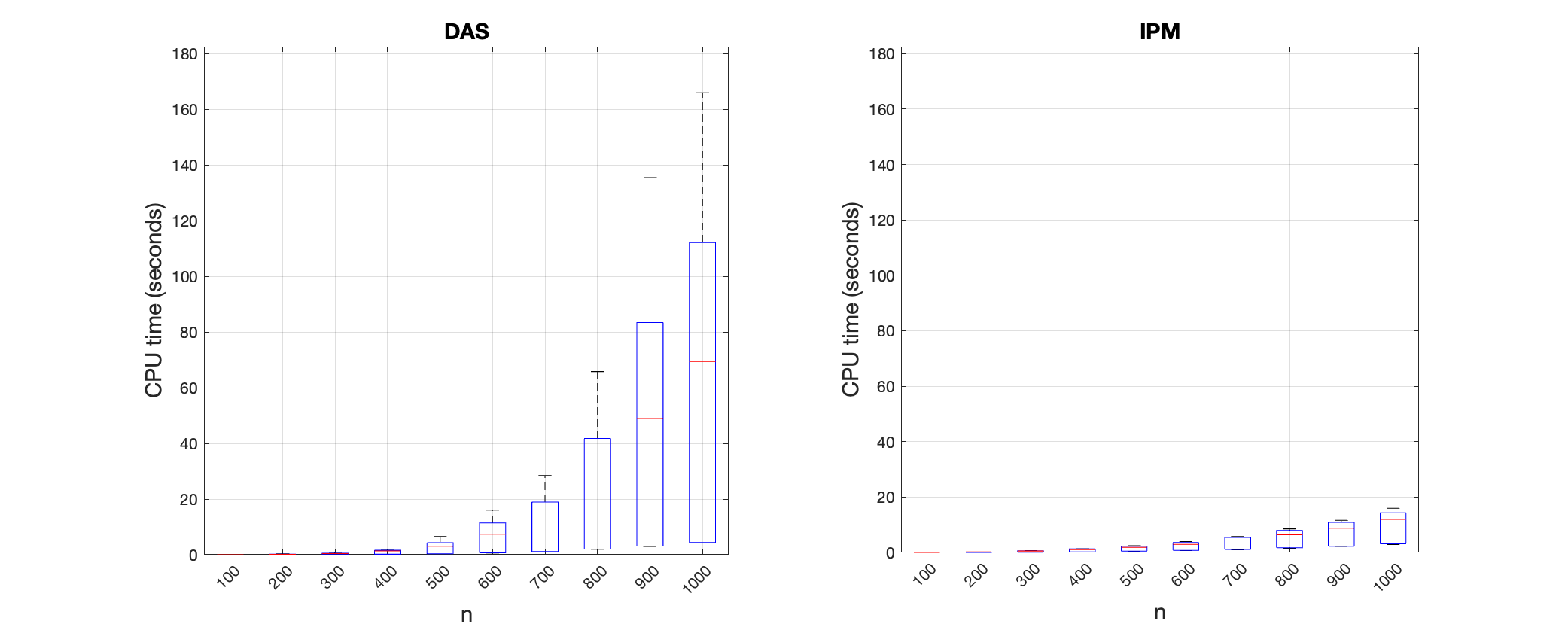}
  \caption{DAS vs. IPM, $m = 2n$}
  \label{fig.das_ipm_6}
\end{figure}

Figures~\ref{fig.das_ipm_4}--\ref{fig.das_ipm_6} parse the results based on different $m$ values (i.e., the number of columns of $G$), namely, $m \in \{n+1,1.5n,2n\}$.  Here, one finds for $n = m+1$ that the DAS solution times are often competitive with the IPM solution times.  However, for DAS there are outliers for which the solution times are much worse, whereas the solution times for IPM are consistently low, even as $n$ increases.  For $m \in \{1.5n,2n\}$, the solution times for DAS are regularly much higher.

One setting not covered in these experiments is $m \ll n$, which in fact is often the case for subproblems that arise in a run of the NonOpt software with its default settings.  Our generation strategy for the QPs in this section---which is based on specifying a particular primal subproblem solution $d_*$---could be modified to allow for $m \ll n$, although the resulting problems would be more contrived since they would require more hand-picked values for $\omega_*$ as opposed to the randomly generated values that are allowed when $m \geq n$.  That said, the results in the subsequent sections demonstrate the use of both DAS and IPM within the algorithms implemented in NonOpt.  By default, the direction computation strategies in NonOpt employ DAS when the number of columns of $G_k$ is less than or equal to 25, and otherwise IPM is employed.  A user can change input options to use only DAS or only IPM, if desired.  The threshold of 25 was chosen based on our experiments in the subsequent sections as a value that led to lower overall solution times.

\subsection{Experiments designed for speed}\label{sec.speed}

Our second set of experiments is designed to demonstrate that with a fixed set of input options, the NonOpt software can generate good solutions of a diverse set of test problems in a computationally efficient manner.  Perhaps the most challenging issue when it comes to the implementation of a general-purpose algorithm for locally Lipschitz minimization is selecting effective termination criteria.  Theoretically, an ideal stationarity condition would be to terminate an algorithm in iteration~$k$ if and only if $\mathrm{dist}_{\gengradient f(x_k)}(0) \leq \epsilon$ for some termination tolerance $\epsilon \in \R{}_{>0}$.  However, in general, it is not reasonable to be able to compute the entire set of generalized gradients $\gengradient f(x_k)$.  NonOpt attempts to approximate the minimum-norm element of this set of generalized gradients through the solution $\omega_k$ of the dual subproblem \eqref{eq.dual_subproblem}.  Indeed, recall that $\ones^T\omega_k = 1$ and $\omega_k \geq 0$, which in turn means that $G_k\omega_k$ is a convex combination of gradients evaluated at points in a small neighborhood of $x_k$.  Therefore, $G_k\omega_k$ can be viewed as an approximation of the minimum-norm element of $\gengradient f(x_k)$.  This means that, in an idealized setting, the algorithm would terminate if and only if $\|G_k\omega_k\|$ is less than a prescribed tolerance.

Practically, however, such an idealized criterion is not reasonable if one wants the software to terminate without excessive computation.  After all, to represent the minimum-norm element in the convex hull of a set in $\R{n}$, one may require (by Caratheodory's theorem) up to $n+1$ gradients.  This means that to enforce such a termination criterion, the algorithm may need to generate $G_k$ with at least $n+1$ columns.  Then, computing the minimum-norm element of the convex hull of the columns of $G_k$ would require a very expensive QP subproblem solve, even when~$n$ is only moderately large.  Therefore, in order to terminate more quickly---at the possible risk of termination before a more accurate solution could be obtained---NonOpt employs termination conditions related to improvement in the objective function value over recent iterations.  (Other software packages for locally Lipschitz minimization also rely heavily on such heuristic termination conditions, including LMBM \cite{Karm2007} and Hanso \cite{Over2010}; see \cite{CurtQue15} for data related to this observation.)  Specifically, NonOpt maintains a counter of how many of the latest iterations resulted in an objective function decrease that was below a prescribed threshold, call it $\Delta f \in \R{}_{>0}$.  If this counter exceeds a prescribed threshold $n_f \in \N{}$ (or if $\|G_k\omega_k\|$ is below a threshold proportional to $\epsilon_k$), then the algorithm sets $\epsilon_{k+1}$ to a fraction of $\epsilon_k$ in line~\ref{step.epsilon} of Algorithm~\ref{alg.nonopt}; otherwise, it sets $\epsilon_{k+1} \gets \epsilon_k$.  Once $\epsilon_k$ is below a prescribed threshold and such a condition holds, then the overall algorithm terminates.  The parameters $\Delta f$ and $n_f$ can be controlled easily by a user.

For our experiments to demonstrate the speed of the NonOpt software with multiple direction computation strategies, we considered the default values for the parameters $\Delta f$ and $n_f$, which are $10^{-5}$ and $10$, respectively.  (The default settings presume that a user would aim for speed when using NonOpt ``out of the box.'')  We considered a set of 20 test problems.  The first ten problems come from \cite{HaarMietMaek04} and the remainder come from \cite{LuksTumaSiskVlceRame02}.  See Table~\ref{tab.problems} for a list of the problem names and their corresponding minimal values of the objective function, where known.  The problems are all scalable in the sense that they are defined for any $n \in \N{}$.  For our experiments, we chose $n=1000$ for all problems.  An initial point $x_0$ is given with each problem definition in \cite{HaarMietMaek04,LuksTumaSiskVlceRame02}.  We used these initial points in our experiments.  Importantly, we note that the first 12 problems are all convex (i.e., all ten problems from \cite{HaarMietMaek04} and the first two from \cite{LuksTumaSiskVlceRame02}) and the remainder are nonconvex.

\begin{table}[ht]
  \centering
  \caption{Problem names, sizes, and minimal objective function values (\texttt{$f_*$}), where known.}
  \label{tab.problems}
  \tiny
  \texttt{
  \btabular{|lcc|}
  \hline
  \multicolumn{1}{|c}{name} & $n$ & $f_*$ \\
  \hline
  \hline
  ActiveFaces        & 1000 & 0 \\
  BrownFunction\_2   & 1000 & 0 \\
  ChainedCB3\_1      & 1000 & $2(n-1) = 1998$ \\
  ChainedCB3\_2      & 1000 & $2(n-1) = 1998$ \\
  ChainedCrescent\_1 & 1000 & 0 \\
  ChainedCrescent\_2 & 1000 & 0 \\
  ChainedLQ          & 1000 & $-\sqrt{2}(n-1) \approx 1412.799$ \\
  ChainedMifflin\_2  & 1000 & $-706.55$ \\
  MaxQ               & 1000 & 0 \\
  MxHilb             & 1000 & 0 \\
  Test29\_2          & 1000 & 0 \\
  Test29\_5          & 1000 & 0 \\
  Test29\_6          & 1000 & - \\
  Test29\_11         & 1000 & - \\
  Test29\_13         & 1000 & - \\
  Test29\_17         & 1000 & - \\
  Test29\_19         & 1000 & - \\
  Test29\_20         & 1000 & - \\
  Test29\_22         & 1000 & - \\
  Test29\_24         & 1000 & - \\
  \hline
  \etabular
  }
\end{table}

The results obtained from NonOpt using its three available direction computation strategies (see page~\pageref{part.dc}) are presented in Table~\ref{tab.speed}.  The table shows the number of iterations, number of function evaluations, number of gradient evaluations, final objective value, and CPU time required (in seconds) for each run.

\begin{table}[ht]
  \centering
  \caption{Results when running NonOpt \emph{for speed} when solving the test problems summarized in Table~\ref{tab.problems}.  Below, \texttt{name} is the problem name, \texttt{dir} is the direction computation strategy employed (\texttt{CP} for cutting plane, \texttt{GC} for gradient combination, and \texttt{G} for gradient, as explained on page~\pageref{part.dc}), \texttt{iters} is the number of iterations required, \texttt{funcs} is the number of function evaluations required, \texttt{grads} is the number of gradient evaluations required, \texttt{f} is the objective function value of the final iterate, and \texttt{CPU} is the time required (in seconds).}
  \label{tab.speed}
  \tiny
  \texttt{
  \btabular{|llrrrrr|}
    \hline
    \multicolumn{1}{|c}{name} & \multicolumn{1}{c}{alg} & \multicolumn{1}{c}{iters} & \multicolumn{1}{c}{funcs} & \multicolumn{1}{c}{grads} & \multicolumn{1}{c}{$f$} & \multicolumn{1}{c|}{CPU (s)} \\
    \hline
    \hline
ActiveFaces & CP & 22 & 414 & 347 & +3.008077e-09 & 0.107499 \\
 & GC & 22 & 414 & 347 & +3.030957e-09 & 0.084499 \\
 & G  & 34 & 458 & 359 & +2.114735e-10 & 0.084612 \\
\hline
BrownFunction\_2 & CP & 107 & 643 & 178 & +3.995409e-08 & 2.327747 \\
 & GC & 107 & 700 & 434 & +1.077945e-04 & 1.978943 \\
 & G  & 135 & 681 & 142 & +1.751286e-06 & 0.398779 \\
\hline
ChainedCB3\_1 & CP & 382 & 2420 & 383 & +2.001144e+03 & 1.862424 \\
 & GC & 370 & 2304 & 371 & +2.001045e+03 & 1.888899 \\
 & G  & 385 & 1689 & 386 & +2.000694e+03 & 1.322669 \\
\hline
ChainedCB3\_2 & CP & 73 & 343 & 83 & +1.998000e+03 & 0.289017 \\
 & GC & 73 & 347 & 83 & +1.998000e+03 & 0.283544 \\
 & G  & 73 & 262 & 80 & +1.998000e+03 & 0.239194 \\
\hline
ChainedCrescent\_1 & CP & 41 & 208 & 53 & +3.250578e-09 & 0.182489 \\
 & GC & 42 & 230 & 58 & +3.773635e-09 & 0.202950 \\
 & G  & 69 & 320 & 81 & +2.846358e-08 & 0.184924 \\
\hline
ChainedCrescent\_2 & CP & 109 & 860 & 169 & +1.452173e-05 & 1.703804 \\
 & GC & 94 & 699 & 252 & +1.075416e-02 & 1.305987 \\
 & G  & 91 & 540 & 98 & +8.126163e-03 & 0.255329 \\
\hline
ChainedLQ & CP & 88 & 499 & 90 & -1.412669e+03 & 0.297360 \\
 & GC & 88 & 498 & 90 & -1.412668e+03 & 0.298095 \\
 & G  & 88 & 391 & 90 & -1.412672e+03 & 0.285419 \\
\hline
ChainedMifflin\_2 & CP & 82 & 646 & 98 & -7.063199e+02 & 0.459288 \\
 & GC & 80 & 677 & 398 & -7.063075e+02 & 1.620687 \\
 & G  & 80 & 441 & 89 & -7.063076e+02 & 0.257966 \\
\hline
MaxQ & CP & 1797 & 5431 & 1889 & +1.756840e-02 & 5.637300 \\
 & GC & 1797 & 5431 & 1889 & +1.756840e-02 & 5.524017 \\
 & G  & 1797 & 3627 & 1889 & +1.756840e-02 & 5.474417 \\
\hline
MxHilb & CP & 198 & 858 & 275 & +2.114027e-04 & 4.502786 \\
 & GC & 198 & 858 & 275 & +2.114027e-04 & 4.487898 \\
 & G  & 198 & 659 & 275 & +2.114027e-04 & 3.796419 \\
\hline
Test29\_2 & CP & 1443 & 7871 & 1555 & +2.024142e-03 & 5.872923 \\
 & GC & 1369 & 7751 & 1580 & +2.744651e-03 & 4.633045 \\
 & G  & 1486 & 6676 & 1557 & +1.803336e-03 & 4.489679 \\
\hline
Test29\_5 & CP & 116 & 332 & 117 & +5.457738e-06 & 2.212484 \\
 & GC & 116 & 332 & 117 & +5.457738e-06 & 2.207088 \\
 & G  & 116 & 215 & 117 & +5.457738e-06 & 1.802998 \\
\hline
Test29\_6 & CP & 111 & 910 & 164 & +2.000000e+00 & 2.774004 \\
 & GC & 102 & 832 & 292 & +2.000211e+00 & 1.722507 \\
 & G  & 125 & 643 & 126 & +2.000000e+00 & 0.384783 \\
\hline
Test29\_11 & CP & 93 & 1110 & 134 & +1.203082e+04 & 0.617999 \\
 & GC & 61 & 598 & 86 & +1.205794e+04 & 0.513414 \\
 & G  & 60 & 475 & 91 & +1.207587e+04 & 0.119095 \\
\hline
Test29\_13 & CP & 1082 & 6813 & 1087 & +5.671352e+02 & 13.372474 \\
 & GC & 79 & 740 & 183 & +5.799427e+02 & 1.332347 \\
 & G  & 161 & 960 & 170 & +5.799165e+02 & 0.945128 \\
\hline
Test29\_17 & CP & 44 & 656 & 81 & +1.091000e-03 & 0.559470 \\
 & GC & 38 & 407 & 129 & +1.092484e-03 & 0.529614 \\
 & G  & 54 & 656 & 66 & +1.092500e-03 & 0.188093 \\
\hline
Test29\_19 & CP & 111 & 954 & 165 & +1.000000e+00 & 3.042539 \\
 & GC & 96 & 835 & 276 & +1.000135e+00 & 1.859920 \\
 & G  & 118 & 612 & 120 & +1.000000e+00 & 0.368168 \\
\hline
Test29\_20 & CP & 87 & 761 & 173 & +5.000018e-01 & 1.393339 \\
 & GC & 92 & 856 & 414 & +5.000235e-01 & 1.940743 \\
 & G  & 128 & 861 & 250 & +5.000000e-01 & 0.402968 \\
\hline
Test29\_22 & CP & 74 & 678 & 123 & +4.531179e-05 & 0.479912 \\
 & GC & 39 & 646 & 235 & +4.533281e-05 & 0.220620 \\
 & G  & 79 & 956 & 170 & +2.190039e-05 & 0.271480 \\
\hline
Test29\_24 & CP & 103 & 920 & 190 & +1.099124e-01 & 2.947687 \\
 & GC & 91 & 849 & 418 & +1.099203e-01 & 1.933301 \\
 & G  & 102 & 619 & 118 & +1.099124e-01 & 0.275045 \\
\hline
  \etabular
  }
\end{table}

A couple of observations are worthwhile to make at this point.  (For further observations, see the next subsection.)  First, observe that, for the most part, the final objective values obtained using the three different direction computation strategies were relatively consistent. Second, observe that the final objective values obtained when solving the convex problems are typically close to their known optimal values; see Table~\ref{tab.problems}.  (For some problems, one might desire additional digits of accuracy, which motivates the experiments in the next subsection.)  Third, observe that the time required for each run was at most a few seconds.  The total time required to run this entire experiment---i.e., all runs of all problems---was approximately 85 seconds.  We comment on this further in the next subsection.

\subsection{Experiments designed for accuracy}\label{sec.accuracy}

We ran a second set of experiments using the same setup as in the prior subsection, but with tighter termination tolerances.  There are various input parameters available to a user of the software that the user can adjust in order to try to obtain a more accurate solution estimate.  For our purposes here, the only adjustment that we made was to set $\Delta f \gets 10^{-8}$ and $n_f \gets 20$ (as opposed to $\Delta f \gets 10^{-5}$ and $n_f \gets 10$ in the previous subsection).  The results are provided in Table~\ref{tab.accuracy}.

\begin{table}[ht]
  \centering
  \caption{Results when running NonOpt \emph{for accuracy} when solving the test problems summarized in Table~\ref{tab.problems}.  Below, \texttt{name} is the problem name, \texttt{dir} is the direction computation strategy employed (\texttt{CP} for cutting plane, \texttt{GC} for gradient combination, and \texttt{G} for gradient, as explained on page~\pageref{part.dc}), \texttt{iters} is the number of iterations required, \texttt{funcs} is the number of function evaluations required, \texttt{grads} is the number of gradient evaluations required, \texttt{f} is the objective function value of the final iterate, and \texttt{CPU} is the time required (in seconds).}
  \label{tab.accuracy}
  \tiny
  \texttt{
  \btabular{|llrrrrr|}
    \hline
    \multicolumn{1}{|c}{name} & \multicolumn{1}{c}{dir} & \multicolumn{1}{c}{iters} & \multicolumn{1}{c}{funcs} & \multicolumn{1}{c}{grads} & \multicolumn{1}{c}{$f$} & \multicolumn{1}{c|}{CPU (s)} \\
    \hline
    \hline
ActiveFaces & CP & 22 & 414 & 347 & +3.008077e-09 & 0.083728 \\
 & GC & 22 & 414 & 347 & +3.030957e-09 & 0.078691 \\
 & G  & 62 & 582 & 387 & +1.156546e-10 & 0.130724 \\
\hline
BrownFunction\_2 & CP & 148 & 1256 & 279 & +2.280562e-08 & 4.199453 \\
 & GC & 493 & 3920 & 5077 & +1.040409e-04 & 23.908808 \\
 & G  & 200 & 991 & 207 & +1.894216e-06 & 0.508164 \\
\hline
ChainedCB3\_1 & CP & 1713 & 11832 & 2252 & +1.998001e+03 & 33.155979 \\
 & GC & 1550 & 10283 & 3452 & +1.998002e+03 & 28.769716 \\
 & G  & 730 & 3530 & 741 & +1.998489e+03 & 2.446144 \\
\hline
ChainedCB3\_2 & CP & 92 & 487 & 120 & +1.998000e+03 & 0.554696 \\
 & GC & 81 & 441 & 112 & +1.998000e+03 & 0.430946 \\
 & G  & 135 & 578 & 143 & +1.998000e+03 & 0.370261 \\
\hline
ChainedCrescent\_1 & CP & 41 & 208 & 53 & +3.250578e-09 & 0.187423 \\
 & GC & 42 & 230 & 58 & +3.773635e-09 & 0.213865 \\
 & G  & 116 & 551 & 142 & +2.849132e-08 & 0.283621 \\
\hline
ChainedCrescent\_2 & CP & 342 & 3460 & 631 & +9.485031e-06 & 12.051088 \\
 & GC & 1177 & 10165 & 12800 & +9.118619e-03 & 68.642038 \\
 & G  & 131 & 712 & 138 & +8.126163e-03 & 0.329044 \\
\hline
ChainedLQ & CP & 751 & 5553 & 1182 & -1.412799e+03 & 21.581284 \\
 & GC & 1233 & 8432 & 1790 & -1.412799e+03 & 24.280076 \\
 & G  & 793 & 3434 & 796 & -1.412799e+03 & 2.595213 \\
\hline
ChainedMifflin\_2 & CP & 2084 & 13223 & 2190 & -7.065434e+02 & 27.780682 \\
 & GC & 133 & 1096 & 944 & -7.063075e+02 & 4.219681 \\
 & G  & 138 & 769 & 151 & -7.063076e+02 & 0.390109 \\
\hline
MaxQ & CP & 1888 & 5825 & 1996 & +9.846592e-03 & 5.841492 \\
 & GC & 1888 & 5825 & 1996 & +9.846592e-03 & 5.823572 \\
 & G  & 1888 & 3930 & 1996 & +9.846592e-03 & 5.896577 \\
\hline
MxHilb & CP & 469 & 1963 & 642 & +2.043153e-06 & 10.497335 \\
 & GC & 469 & 1963 & 642 & +2.043153e-06 & 10.470230 \\
 & G  & 469 & 1493 & 642 & +2.043153e-06 & 8.779602 \\
\hline
Test29\_2 & CP & 2899 & 17619 & 3471 & +1.253118e-06 & 21.473363 \\
 & GC & 2440 & 16181 & 3217 & +1.861335e-05 & 14.838679 \\
 & G  & 2948 & 14754 & 3119 & +1.176407e-06 & 9.079778 \\
\hline
Test29\_5 & CP & 202 & 732 & 203 & +8.370109e-07 & 4.356180 \\
 & GC & 202 & 732 & 203 & +8.370109e-07 & 4.350333 \\
 & G  & 202 & 529 & 203 & +8.370109e-07 & 3.661412 \\
\hline
Test29\_6 & CP & 157 & 1301 & 255 & +2.000000e+00 & 4.675574 \\
 & GC & 2458 & 17606 & 5311 & +6.923490e-02 & 19.347360 \\
 & G  & 231 & 1156 & 232 & +2.000000e+00 & 0.612321 \\
\hline
Test29\_11 & CP & 205 & 1937 & 295 & +1.203129e+04 & 4.611533 \\
 & GC & 1436 & 10319 & 7261 & +1.202995e+04 & 40.581987 \\
 & G  & 117 & 760 & 121 & +1.207587e+04 & 0.214988 \\
\hline
Test29\_13 & CP & 2102 & 12392 & 2141 & +5.661129e+02 & 23.540974 \\
 & GC & 1723 & 10409 & 2385 & +5.661128e+02 & 22.928162 \\
 & G  & 1962 & 8684 & 1971 & +5.661131e+02 & 10.505420 \\
\hline
Test29\_17 & CP & 262 & 4153 & 504 & +1.057500e-03 & 5.205452 \\
 & GC & 1196 & 13746 & 3294 & +3.668662e-05 & 23.924566 \\
 & G  & 179 & 1902 & 215 & +1.080500e-03 & 0.619420 \\
\hline
Test29\_19 & CP & 149 & 1279 & 241 & +1.000000e+00 & 4.645676 \\
 & GC & 1033 & 8520 & 4538 & +2.543375e-01 & 20.473368 \\
 & G  & 214 & 1077 & 216 & +1.000000e+00 & 0.576925 \\
\hline
Test29\_20 & CP & 215 & 2088 & 584 & +5.000000e-01 & 8.460440 \\
 & GC & 622 & 7547 & 7004 & +5.000085e-01 & 38.719123 \\
 & G  & 284 & 1742 & 524 & +5.000000e-01 & 0.832635 \\
\hline
Test29\_22 & CP & 421 & 5124 & 1232 & +8.267981e-06 & 12.493644 \\
 & GC & 216 & 3511 & 1008 & +7.592961e-06 & 3.728066 \\
 & G  & 574 & 6037 & 1289 & +1.501467e-06 & 1.941063 \\
\hline
Test29\_24 & CP & 145 & 1271 & 273 & +1.099124e-01 & 4.663748 \\
 & GC & 468 & 5633 & 5036 & +1.099154e-01 & 23.860827 \\
 & G  & 166 & 910 & 182 & +1.099124e-01 & 0.368457 \\
\hline
  \etabular
  }
\end{table}

A few observations are worthwhile to make, especially when one compares the results in Table~\ref{tab.accuracy} to those in Table~\ref{tab.speed}.  First, it should be said that for some problems the gains in terms of final objective values are only modest despite the tighter tolerances and (in some cases significantly) increased computational effort.  This suggests that the default settings of NonOpt---which are designed for speed---can typically yield good final objective values.  Second, observe that in other cases the final solution quality improves somewhat significantly with the tighter tolerances.  Third, observe that in some nonconvex instances the direction computation strategy \texttt{GC} appears to find a solution that is substantially better than even those obtained by the other strategies.  With respect to this, we draw the reader's attention to the results for \texttt{Test29\_6}, \texttt{Test29\_17}, and \texttt{Test29\_19}.  For these problems, strategy \texttt{GC} with the tighter tolerances led to final objective values that are improved by around an order of magnitude or more.  All of the solution quality improvements in these experiments come at a cost, however.  Whereas the total time required for the experiments in the previous subsection was around 85 seconds, the total time required for the experiments here was around 580 seconds.

\subsection{Experiments on large-scale image-processing problems}

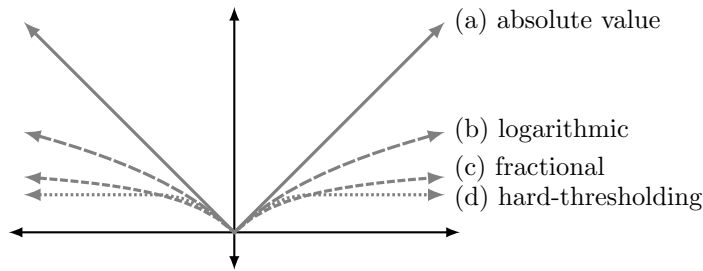
\begin{figure}[ht]
  \centering
  \begin{tikzpicture}
  \coordinate (xmax) at ( 3.0, 0.0);
  \coordinate (ymax) at ( 0.0, 3.0);
  \coordinate (xmin) at (-3.0, 0.0);
  \coordinate (ymin) at ( 0.0,-0.5);
  \draw[black,thick,latex-latex] (xmin) -- (xmax);
  \draw[black,thick,latex-latex] (ymin) -- (ymax);
  \draw[gray,very thick,latex-latex,domain=-2.8:2.8                            ] plot({\x},{1.0*abs(\x)});
  \draw[gray,very thick,latex-latex,domain=-2.8:2.8,dash pattern=on 5pt off 1pt] plot({\x},{ln(1 + 1.0*abs(\x))});
  \draw[gray,very thick,latex-latex,domain=-2.8:2.8,dash pattern=on 3pt off 1pt] plot({\x},{1.0*abs(\x)/(1 + 1.0*abs(\x))});
  \draw[gray,very thick,latex-latex,domain=-2.8:2.8,dash pattern=on 1pt off 1pt] plot({\x},{1.0/2.0 - max(0,1.0 - abs(\x))*max(0,1.0 - abs(\x))/(2.0*1.0)});
  \coordinate [label=right:(a) absolute value] (a) at (2.8,2.8);
  \coordinate [label=right:(b) logarithmic] (l) at (2.8,1.35);
  \coordinate [label=right:(c) fractional] (f) at (2.8,0.85);
  \coordinate [label=right:(d) hard-thresholding] (h) at (2.8,0.45);
  \end{tikzpicture}
  \caption{Regularizers for image denoising problem instances, each illustrated for $\beta = 1$.}
  \label{fig.regularizers}
\end{figure}

Our final set of experiments is designed to demonstrate that NonOpt is able to solve large-scale problems efficiently, and that it can transition seamlessly from solving convex to nonconvex problems.  For these experiments, we considered a set of related image-denoising problems.  Suppose that a grayscale image represented by a matrix $X_o \in [256]^{n_r} \times [256]^{n_r}$ has been corrupted by noise so that one only has access to a noisy grayscale image $\Xbar \in [256]^{n_r} \times [256]^{n_r}$ with $\Xbar \approx X_o$.  (The presence of $[256]$ in these spaces is due to the fact that each pixel value is represented by an integer between 1 and 256.  Moreover, $n_r \in \N{}$ is the number of rows in the image and $n_c \in \N{}$ is the number of columns.)  The goal is to denoise the image by solving an optimization problem.  To demonstrate such a procedure, the optimization problems that we consider for our experiments each have the form
\bequation\label{eq.denoising}
  \baligned
    \min_{X \in \R{n_r \times n_c}} &\ \sum_{i\in[n_r]} \sum_{j\in[n_c]} (X_{i,j} - \Xbar_{i,j})^2 \\
    &\ + \lambda \sum_{i \in [n_r-1]} \sum_{j \in [n_c-1]} \phi_\beta(X_{i+1,j} - X_{i,j}) + \phi_\beta(X_{i,j+1} - X_{i,j}),
  \ealigned
\eequation
where $\lambda \in \R{}_{>0}$ is a regularization parameter and $\phi_\beta : \R{} \to \R{}$ is a regularization function (also known as a regularizer) that is parameterized by $\beta \in \R{}_{>0}$.  (As is common in the literature on image denoising, problem~\eqref{eq.denoising} is defined over $\R{n_r \times n_c}$, although to render the solution as an image we round each pixel value to the nearest integer in $[256]$.)  Overall, we consider the following four choices of the regularizer:
\benumerate
  \item[(a)] absolute value: $\phi_\beta(t) = \beta|t|$;
  \item[(b)] logarithmic: $\phi_\beta(t) = \beta\log(1 + \beta |t|)$;
  \item[(c)] fractional: $\phi_\beta(t) = \frac{\beta |t|}{1 + \beta |t|}$; and
  \item[(d)] hard thresholding: $\phi_\beta(t) = \frac{\beta}{2} - \frac{\max \{ 0, \beta - |t|\}^2}{2\beta}$. \label{en.regularizers}
\eenumerate
The absolute value regularizer can be derived as the anisotropic version of the well known total variation regularizer \cite{RudiOsheFate1992}.  The use of this particular regularizer makes the objective of problem~\eqref{eq.denoising} convex.  Each of the other three choices of the regularizer makes the objective of problem~\eqref{eq.denoising} nonconvex; see Figure~\ref{fig.regularizers}.

We remark in passing that there exist numerous specialized algorithms for solving image denoising problems.  The aim of our experiments here is not to attempt to compete with the state-of-the-art methods for solving such problems.  Rather, our aim is to demonstrate that NonOpt can be employed to solve large-scale problem instances and to highlight the fact that NonOpt can be employed to solve problem~\eqref{eq.denoising} regardless of the choice of regularizer---convex or nonconvex.

\begin{figure}[ht]
  \centering
  \includegraphics[width=0.48\textwidth]{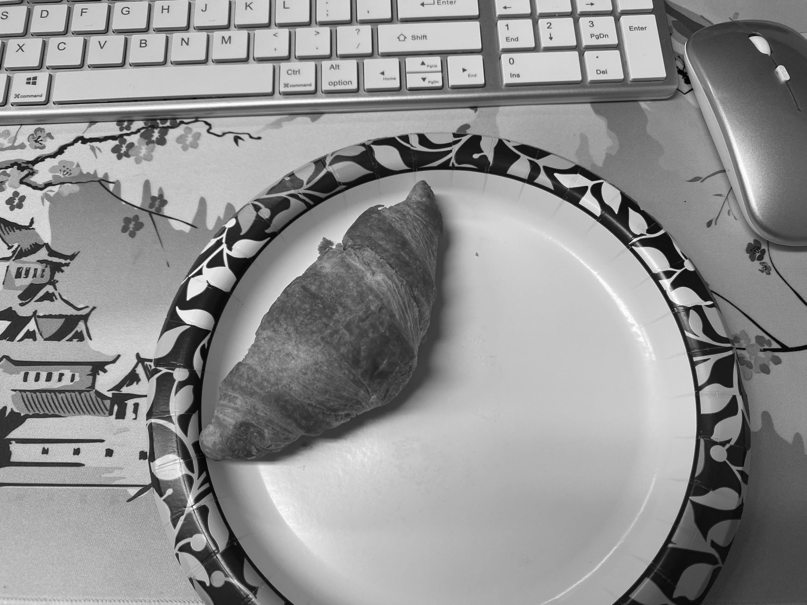}\ 
  \includegraphics[width=0.48\textwidth]{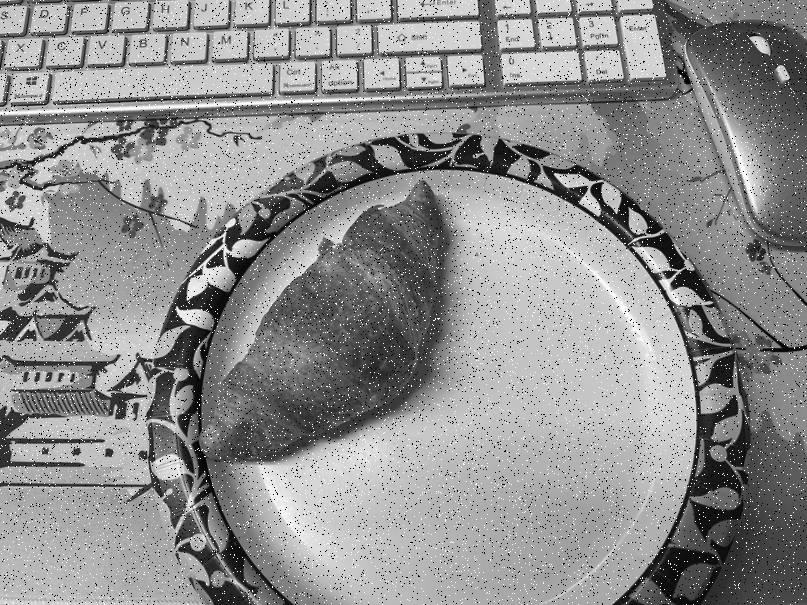}
  \caption{Original (left) and noisy (right) image for denoising experiments.}
  \label{fig.originals}
\end{figure}

We considered $X_o$ as a 605-by-807 grayscale image, the original of which is shown on the left-hand side of Figure~\ref{fig.originals}.  We generated $\Xbar$ by adding ``salt and pepper'' noise to $X_o$ using the default settings of \texttt{imnoise} in Matlab's image processing toolbox.  The resulting noisy image is shown on the right-hand side of Figure~\ref{fig.originals}.

\begin{figure}[ht]
  \centering
  \includegraphics[width=0.48\textwidth]{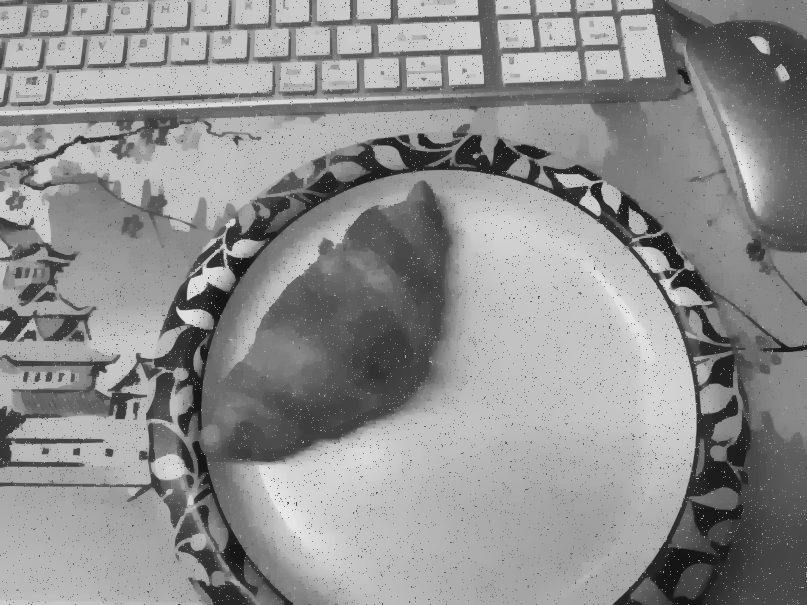}
  \includegraphics[width=0.48\textwidth]{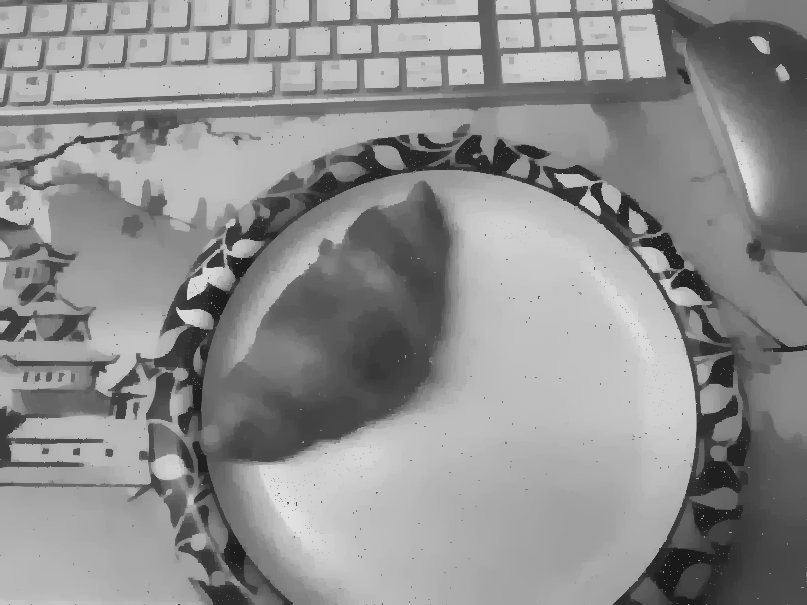} \\
  \includegraphics[width=0.48\textwidth]{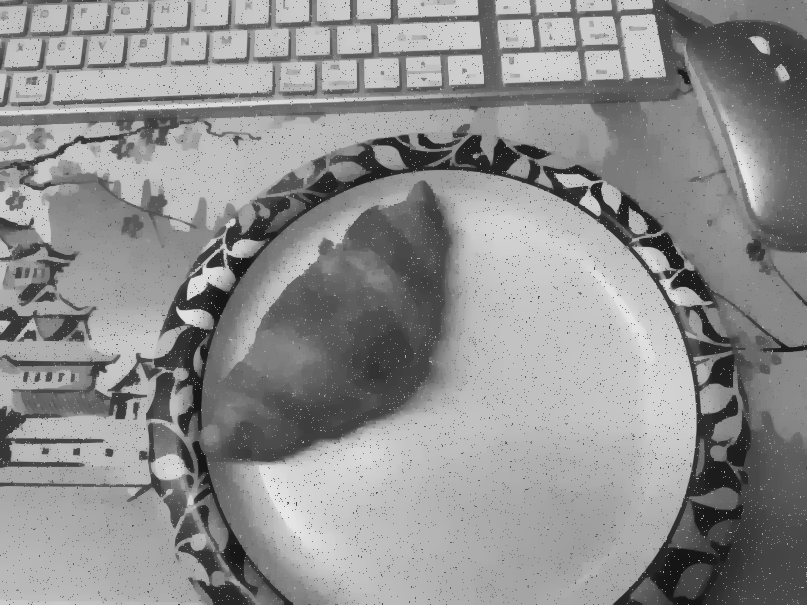}
  \includegraphics[width=0.48\textwidth]{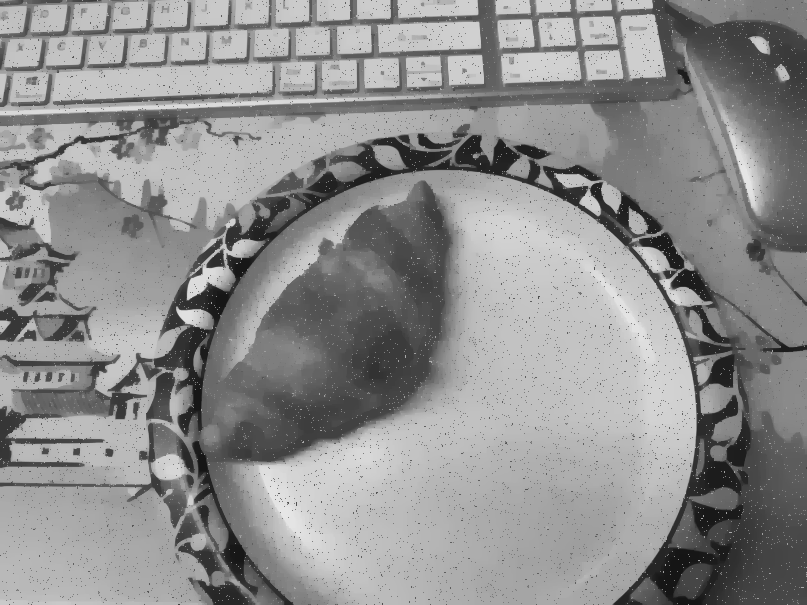}
  \caption{Recovered images obtained with regularizer (a) absolute value (top, left), (b) logarithmic (top, right), (c) fractional (bottom, left), and (d) hard thresholding (bottom, right).}
  \label{fig.recovered}
\end{figure}

We employed NonOpt to solve instances of problem~\eqref{eq.denoising} using each of the four aforementioned regularizers.  Given the size of the image that we considered, each instance had 488,235 variables.  Due to the large scale of these problems, we used the limited memory \texttt{SymmetricMatrix} strategy, i.e., L-BFGS for the Hessian approximations; see page~\pageref{part.sm}.  In order to present results for reasonable scenarios, we tuned $(\lambda,\beta) \in \R{}_{>0} \times \R{}_{>0}$ for each instance by searching over pairs of parameter values of the form $(2^i,2^j)$ for all $(i,j) \in [-20,25] \times [-20,25]$ and, for each instance, choosing the values that led to the smallest mean-squared error (MSE) with the original image $X_o$.  The resulting tuned input parameter values, MSE values, and solution time required for each instance are shown in Table~\ref{tab.mse}.  (For each instance, we take the final solution obtained by NonOpt, round each value to the nearest integer in $[256]$ to produce a matrix $X_*$, and compute the resulting MSE with $X_o$.)  The resulting denoised images are shown in Figure~\ref{fig.recovered}.  Overall, given the large scale of the problem instances, NonOpt is able to find good-quality solutions to both convex and nonconvex instances in relatively little time.

\btable
  \centering
  \caption{Data for solving instances of problem~\eqref{eq.denoising}.  Below, \texttt{regularizer} is the regularizer employed (see page~\pageref{en.regularizers}), $\lambda$ is the regularization parameter, $\beta$ is the regularizer parameter, \texttt{MSE} is the final mean squared error of the solution obtained $X_*$ with respect to the true image $X_o$, and \texttt{CPU} is the number of seconds required to compute $X_*$.}
  \label{tab.mse}
  \texttt{
  \btabular{ccccc}
    \hline
    regularizer & $\lambda$ & $\beta$ & MSE: $\|X_* - X_o\|_F^2$ & CPU (s)\\
    \hline
    (a) & $2^6$    & $2^0$    & 1.08781e+08 & 145.271 \\
    (b) & $2^{21}$ & $2^{-7}$ & 1.20072e+08 & 156.497 \\
    (c) & $2^{25}$ & $2^{-19}$ & 1.08780e+08 & 110.409 \\
    (d) & $2^6$    & $2^{18}$ & 1.08788e+08 & 134.867 \\
    \hline
  \etabular
  }
\etable

\section{Conclusion}\label{sec.conclusion}

We have provided a thorough introduction to an open-source C++ software package for solving locally Lipschitz minimization problems.  We have discussed the overall algorithmic framework that is implemented in the software, which includes both gradient-sampling and proximal-bundle methodologies, and provided an overview of the various modules in the software that allow a user to tailor the solver to their particular application(s).  We have also discussed the details of a newly implemented interior-point subproblem solver that leads to significantly lower solution times for each subproblem and the overall NonOpt algorithmic strategy in general.  The results of our numerical experiments have shown the benefits of the new QP subproblem solver, the speed and reliability of the software, and that NonOpt can solve large-scale (convex and nonconvex) problems efficiently.

\bibliographystyle{plain}
\bibliography{references}

\end{document}